\documentclass{article}

\usepackage[utf8]{inputenc}
\usepackage{amsmath,amssymb,amsthm,bm}
\usepackage{amsfonts}
\usepackage{mathrsfs}
\usepackage{graphicx}
\usepackage{epstopdf}
\usepackage{enumitem}
\usepackage{overpic}
\usepackage{multirow}
\usepackage{nicefrac}
\numberwithin{equation}{section}
\usepackage{float}
\usepackage{caption}
\usepackage{multicol}
\usepackage{dsfont}
\usepackage{bbm}
\usepackage{url}
\usepackage{verbatim}
\usepackage{algorithm,algorithmic}

\usepackage{xcolor}
\usepackage{arydshln}
\usepackage{subfig}
\usepackage{nicematrix}
\usepackage{pifont}

\usepackage[colorlinks,linkcolor=blue,
            citecolor=blue]{hyperref}

\usepackage{geometry}
\newcommand{\trace}{\mathrm{trace}}

\newcommand{\vectx}{{\text{vec}}}

\newtheorem{theorem}{Theorem}

\newtheorem{lemma}{Lemma}

\newtheorem{remark}{Remark}

\providecommand{\keywords}[1]
{
  \small	
  \textbf{\textit{Key words---}} #1
}

\graphicspath{
{images/}
}

\title{A Novel Double Periodic Conformal Flattening Algorithm for Genus-One Surfaces}
\author{Zhong-Heng Tan\thanks{Department of Mathematics, The Chinese University of Hong Kong, Hong Kong. Email: {zhtan@math.cuhk.edu.hk} }
\and Tiexiang Li\thanks{Corresponding Author. School of Mathematics and Shing-Tung Yau Center, Southeast University,
Nanjing 210096; Shanghai Institute for Mathematics and Interdisciplinary Sciences (SIMIS), Shanghai 200433, People's Republic of China. Email: txli@seu.edu.cn}
\and Wen-Wei Lin\thanks{Shanghai Institute for Mathematics and Interdisciplinary Sciences (SIMIS), Shanghai 200433, People’s Republic of China. Email: wwlin@outlook.com}
\and Lok Ming Lui\thanks{Department of Mathematics, The Chinese University of Hong Kong, Hong Kong. Email: lmlui@math.cuhk.edu.hk}
}
\date{}

\begin{document}

\maketitle

\begin{abstract}
    In this paper, we propose novel parameterization methods for genus-one surfaces, called the Double Periodic Conformal Flattening (DPCF) algorithm.
The desired conformal map is obtained by minimizing a conformal energy functional under periodic boundary conditions, which is characterized as an easily solvable quadratic functional minimization problem, yielding a sparse linear system.
The proposed DPCF algorithm offers several key advantages:
(a) the optimal boundary and periodic translation vectors are obtained simultaneously with the conformal map;
(b) the resulting map is independent of the chosen cutting path, thus introducing no extra conformal distortion near the cut seams;
(c) bijectivity is guaranteed under the positive edge weights condition, which can be satisfied by, e.g., using an intrinsic Delaunay triangulation.
Based on this guaranteeing, a simple strategy is employed to ensure bijectivity of the resulting maps for general triangulations.
Numerical experiments illustrate that DPCF algorithm exhibits high accuracy and a 5 fold improvement over the state-of-the-art algorithm in terms of efficiency.
Applications on texture mapping and medical imaging illustrate the practicality of our developed algorithm.
\end{abstract}

\section{Introduction}

The advent of high-precision 3D technology has given rise to a plethora of surfaces that are characterized by high resolutions and complex geometries, which poses challenges to direct surface manipulation methods.
Surface flattening is a fundamental technique that aims to flatten a surface onto a plane, transforming the task of processing a complicated surface into a simple 2D domain-related procedure.
Generally, flattening maps should be bijective or locally injective to ensure optimal visual effects and high processing performance in practical applications, which are the focuses of some works \cite{VGLK14,JPCY20,QWWZ23}.
The common goal of surface flattening is to reduce the degree of flattening distortion as much as possible.
Conformal flattening is a flattening technique for minimizing conformal distortion, which preserves the local shape by the conformal map.
Consequently, this approach plays a significant role in computer-aided engineering \cite{ULKH00,WCSS24} and computer graphics \cite{SHSA00,LBPS02}.

For the genus-zero surface with boundary,
various algorithms and methods have been developed to achieve conformal surface flattening, including most isometric parametrizations (MIPS) \cite{Horm00}, angle-based flattening \cite{SASE01,SALB05}, LSCM \cite{LBPS02}, DCP \cite{MDMM02}, circle patterns \cite{LKBS06}, spectral conformal parameterization \cite{PMYT08}, quasi-conformal approach \cite{PTLM15,GPLM18,Choi21} and discrete conformal equivalence \cite{MGBS21}.
Ricci flow \cite{MJJK08,YLRG09} is an effective method for computing conformal factors and subsequently constructing the conformal map.
Similarly, the Yamabe equation \cite{RSKC17} is also used to flatten surfaces with various boundary conditions.
Conformal energy minimization has been proposed for disk parameterizations \cite{MHWW17,YCWW21} and free boundary flattening \cite{WQXD14-SCP}.

Closed surfaces such as spheres and tori can be cut along a graph tree to be disk-topology, so that the cut surfaces are flattened to $\mathbb{R}^2$ \cite{AMDZ12,NARP15,YSDS18,QFWO21}.
In this scenario, the cut seam correspondence should be taken into consideration.
Generally, each pair pf boundaries from the same seam should have identical shape, i.e., they differ only from a rigid transformation.
These maps are called seamless parameterization.
Harmonic map based method \cite{NAYL15,ABEC17,EHEC19} is a common approach to computing the seamless parameterization; that is, generate weights for each edge and solve a discrete harmonic equation with cut seam constraints.
Floater \cite{Floa02} proved that the harmonic map is bijective if edge weights are positive and the flattened domain is convex.
The positive edge weights could be satisfied by restricting intrinsic Delaunay mesh \cite{NAYL15} or modifying the edge weights, such as mean value coordinates \cite{Floa03}.
However, we can hardly obtain a convex domain under the seam constraints.
For this problem, it is common to add bijective constriants \cite{ABEC17,EHEC19} or removing non-bijective elements during iterating computation \cite{JPXM19,QWWZ23} to obtain locally injective maps or bijective maps.
Furthermore, the boundary condition is not sufficient to guarantee the uniqueness of the solution.
Then additional conditions and constraints, e.g., fixing rigid transformation \cite{NAYL15} and some vertices \cite{ABEC17,EHEC19}, are prescribed, which may lead to larger distortion.

On the other hand,
Gu et al. \cite[Chap. 11]{XGST07book} and Yueh et al. \cite{MHTL20} utilized the de Rham cohomology approach to compute conformal parameterizations for genus-one surfaces, respectively, which computes the associated holomorphic 1-forms first, and then obtains the corresponding holomorphic (conformal) maps by integrating along a suitable path from a fixed starting point to each vertex.
From a theoretical point of view, this approach can be used to understand the relationship between holomorphic 1-forms and maps, and to produce a relatively accurate conformal map on the basis of rigorous theoretical guarantee.
Gortler et al. \cite{SGCG06} proved the bijectivity of holomorphic 1-form methods.
However, from the numerical perspective of finding the conformal map, the computations of holomorphic 1-forms and subsequent line integrals on
the de Rham cohomology theory \cite{MHTL20,XGST07book} is truly a detour.

In this paper, we aim to compute conformal seamless parameterizations for genus-one surfaces via conformal energy minimization and prove the bijectivity of our proposed algorithm.
Since the flattening map produced with optimal translation correspondence for genus-one surfaces leads to a periodic tiling on $\mathbb{R}^2$, we call the proposed approach `double periodic conformal flattening' (DPCF).
The main contributions of this paper are threefold.

\begin{itemize}
    \item We propose an efficient and reliable numerical method, namely, DPCF for genus-one surfaces.
    Based on conformal energy minimization with identical cut seam constraints, the optimal boundary associated with translation vectors are obtained simultaneously with the map by only solving a sparse linear systems with two-column right hand side.
    The proposed algorithms have better performance in terms of accuracy and, in particular, provide efficiency improvements of approximately 5 times over the recently developed state-of-the-art algorithm \cite{MHTL20}.

    \item We show that the maps resulting from the developed algorithm are independent of the selected cutting path and that no additional conformal error is induced near the cutting path, which reflects the stability of our method.
    Based on this result, we further prove that the translation vectors only depend on the homology basis to which the cutting paths belong.

    \item
    We generalize the bijective result in \cite{Floa02} and prove that the computed conformal map is bijective under the positive edge weights condition, e.g., the surface is an intrinsic Delaunay triangulation (IDT).
    We further develop a bijectivity guarantee strategy for non-Delaunay triangulations. Numerical experiments illustrate the significant effectiveness and accuracy improvement of the strategy.
\end{itemize}

This paper is organized as follows.
In Section \ref{sec:surface}, we briefly introduce the concept of a discrete surface and review conformal energy functional.
In Section \ref{sec:DPCF}, we derive conformal energy minimization with periodic boundary conditions and propose DPCF algorithm.
In Section \ref{sec:cut} and \ref{sec:bijective}, we prove the cutting path independence and bijectivity guarantee of the proposed algorithms DPCF, respectively.
The numerical performance and applications of the proposed approaches are presented in Section \ref{sec:NE} and \ref{sec:app}, respectively.
Concluding remarks are given in Section \ref{sec:conclusion}.

The notations that are frequently used in this paper are listed here.
Bold letters, e.g., $\mathbf{f},\mathbf{g}$, denote vectors whose entries are scalar or vertex components.
Furthermore, $[\mathbf{f}]$ denotes the corresponding index set of components of $\mathbf{f}$.
$\mathbf{1}_{m\times n}$ denotes the $m\times n$ matrix of all ones. The same notation without the subscript $\mathbf{1}$ denotes the vector of all ones with appropriate dimensions.
$\mathbf{e}_i$ denotes the $i$-th column of the identity matrix with appropriate dimensions.
$\big[ A \big]_{ij}$ denotes the $(i,j)$-th entry of matrix $A$.
$\widehat{pq}$ denotes a specified simple polyline from vertex $p$ to vertex $q$ along the associated triangular meshes.

\section{Discrete Surfaces and Conformal Energy} \label{sec:surface}

In this section, we briefly introduce discrete surfaces and review the concept of conformal energy.
Let $\mathcal{M}$ be a discrete surface, which consists of sets of vertices
\begin{align*}
    \mathcal{V}(\mathcal{M}) = \{v_\ell = (v_\ell^1, v_\ell^2, v_\ell^3) \in \mathbb{R}^3\}_{\ell=1}^n,
\end{align*}
edges
\begin{align*}
    \mathcal{E}(\mathcal{M}) = \{[v_i,v_j] \subset \mathbb{R}^3~|~ \text{ for some } \{v_i,v_j\} \in \mathcal{V}(\mathcal{M}) \},
\end{align*}
and triangular faces
\begin{align*}
    \mathcal{F}(\mathcal{M}) = \{[v_i,v_j,v_k] \subset \mathbb{R}^3~|~\text{for some } \{v_i, v_j, v_k\} \subset \mathcal{V}(\mathcal{M}) \},
\end{align*}
where the bracketed terms $[v_i,v_j]$ and $[v_i, v_j, v_k]$ denote an edge and a triangle, respectively.
$|[v_i,v_j,v_k]|$ is the area of triangle $[v_i,v_j,v_k]$.

The discrete map $f: \mathcal{M} \to \mathbb{R}^2$ defined on $\mathcal{M}$ is considered a piecewise affine map. From this perspective, the image location $f(v)$ of a point $v$ in the triangle $[v_i,v_j,v_k]
\subset \mathcal{M}$ can be represented by the barycentric coordinates
\begin{align}
    f(v) = \frac{|[v,v_j,v_k]|}{|[v_i,v_j,v_k]|} f(v_i) + \frac{|[v_i,v,v_k]|}{|[v_i,v_j,v_k]|} f(v_j) + \frac{|[v_i,v_j,v]|}{|[v_i,v_j,v_k]|} f(v_k). \label{eq:bary}
\end{align}
Let $\mathbf{f}_{\ell} = f(v_\ell) \in \mathbb{R}^2$ for $\ell = 1,2,\cdots,n$, and define the image vertex matrix
$
\mathbf{f}:= [\mathbf{f}_{1}^\mathrm{T},\mathbf{f}_{2}^\mathrm{T},\cdots,\mathbf{f}_{n}^\mathrm{T}]^\mathrm{T}.
$
The map $f$ is induced by $\mathbf{f}$ via the barycentric coordinates described in \eqref{eq:bary}.

Without confusion, $\mathcal{M}$ can be used for both continuous and discrete cases. The continuous conformal energy functional \cite{Hutc91} for a map $f$ on a 2-manifold $\mathcal{M}$ is defined as
\begin{align*}
    E_C(f) = \frac{1}{2} \int_\mathcal{M} \|\nabla_\mathcal{M} f\|_F^2 ds - |f(\mathcal{M})|,
\end{align*}
where $\nabla_\mathcal{M}$ represents the tangent gradient and $|f(\mathcal{M})|$ is the area of the image surface $f(\mathcal{M})$. The conformal energy satisfies $E_C(f) \geq 0$ for an arbitrary map $f$ and $E_C(f) = 0$ if and only if $f$ is conformal.

The discrete conformal energy \cite[Chap. 16]{XGST07book} for a piecewise affine map $f$ defined on the discrete surface $\mathcal{M}$ is given by
\begin{align}
    E_C(f) = \frac{1}{2} \trace(\mathbf{f}^\mathrm{T} L_D \mathbf{f}) - |f(\mathcal{M})|, \label{eq:EC_discrete}
\end{align}
where $L_D$ is the Laplacian matrix with
\begin{align}
    \big[ L_D \big]_{ij} = \begin{cases}
        -w_{ij} \equiv -\frac{1}{2} (\cot \theta_{ij} + \cot \theta_{ji}), & \text{if } [v_{i},v_j] \in \mathcal{E}(\mathcal{M}),\\
        \sum_{\ell \neq i} w_{i\ell}, & \text{if } i = j, \\
        0, & \text{otherwise},
    \end{cases} \label{eq:LD}
\end{align}
in which $\theta_{ij}$ and $\theta_{ji}$ are the opposite angles of edge $[v_i,v_j]$. Numerically, the discrete conformal map can be computed by minimizing the conformal energy of \eqref{eq:EC_discrete}.

\section{Double Periodic Conformal Flattening for Genus-One Surfaces} \label{sec:DPCF}

Let $\mathcal{M}$ be a genus-one surface. We first illustrate the uniformization theorem \cite{XGST07book} for the classification of a simply connected Riemann surface on a conformal map.

\begin{figure}[htp]
    \centering
    \includegraphics[width=0.5\linewidth]{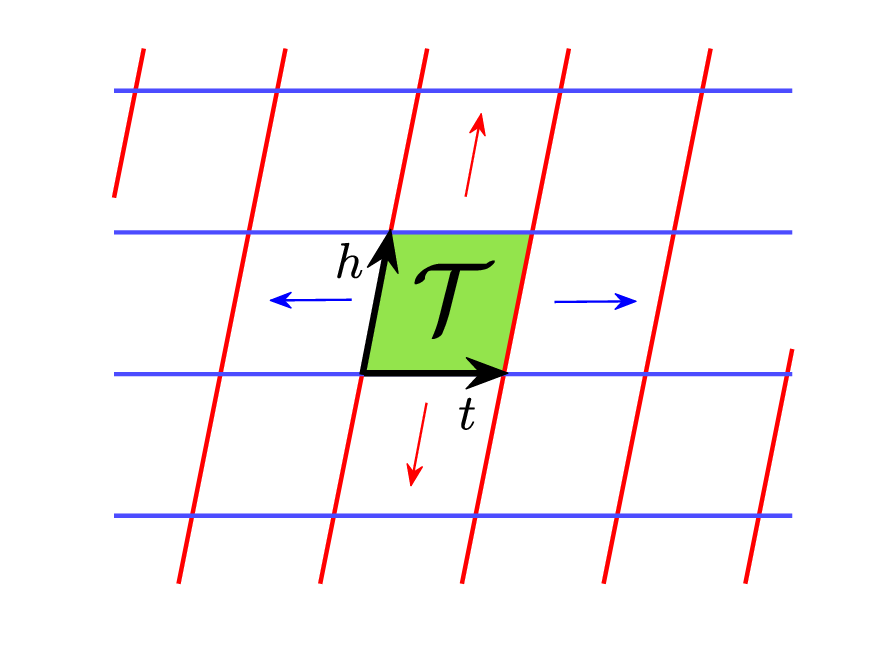}
    \caption{A 2D lattice and the unit cell domain $\mathcal{T}$ with two periodic boundary vectors $h$ and $t$ as in \eqref{eq:lattice2}.}
    \label{fig:lattice}
\end{figure}

\begin{theorem}[Poincar\'e-Klein-Koebe Uniformization Theorem \cite{XGST07book}]
    A simply connected Riemann surface $\mathcal{M}$ is conformally equivalent to one of the following three canonical Riemann surfaces:
    \begin{itemize}
        \item[1.] An extended complex plane $\overline{\mathbb{C}} = \mathbb{C} \cup \{\infty\}$ with $\mathcal{M}$ of genus $=$ $0$;
        \item[2.] A complex plane $\mathbb{C}$ with $\mathcal{M}$ of genus $=$ $1$;
        \item[3.] An open unit disk $\Delta = \{z \in \mathbb{C} | |z| < 1\}$ with $\mathcal{M}$ of genus $>$ $1$.
    \end{itemize}
\end{theorem}
Furthermore, the genus-one surface $\mathcal{M}$ is conformally equivalent to a quotient group $\mathbb{C} / G$, where $G$ is the deck transformation group
\begin{align} \label{group:deckg1}
    G = \{z \to z + k_1\omega_1+k_2\omega_2, k_1,k_2\in \mathbb{Z}, \omega_1,\omega_2 \in \mathbb{C}, \mathrm{Im}(\omega_1/\omega_2) \neq 0\}.
\end{align}
This means that a genus-one surface can be conformally transformed into the unit cell of a 2D lattice
\begin{align}
    \Lambda_{h,t} = \{k_1h+k_2t, k_1,k_2\in \mathbb{Z}, h, t \in \mathbb{R}^{1\times 2}, \det ([h^\mathrm{T}, t^\mathrm{T}]) \neq 0\}. \label{eq:lattice2}
\end{align}
Geometrically, there exists a homeomorphism that transforms $\mathcal{M}$ into a domain $\mathcal{T}$ with two periodic boundary vectors $h$ and $t$, as in \eqref{eq:lattice2}.
This domain $\mathcal{T}$, as a unit cell, can be used to tile the entire 2D plane $\mathbb{R}^2$ via translation along the lattice $\Lambda_{h,t}$, as illustrated in Figure \ref{fig:lattice}, where the blue and red lines correspond to the handle ($h$) and tunnel ($t$) loops, respectively, which are in line with the homology basis of a genus-one surface.
Our goal is to find a conformal map that can transform $\mathcal{M}$ into the unit cell domain $\mathcal{T}$.
Since a unit cell has two pairs of periodic sides, this approach is called double periodic conformal flattening.

\begin{figure}[htp]
    \centering
\subfloat[A genus-one surface with the handle loop $\alpha$ and the tunnel loop $\beta$]{\label{subfig:G1sufHT}
\includegraphics[clip,trim = {3cm 2cm 3cm 2cm},width = 0.3\textwidth]{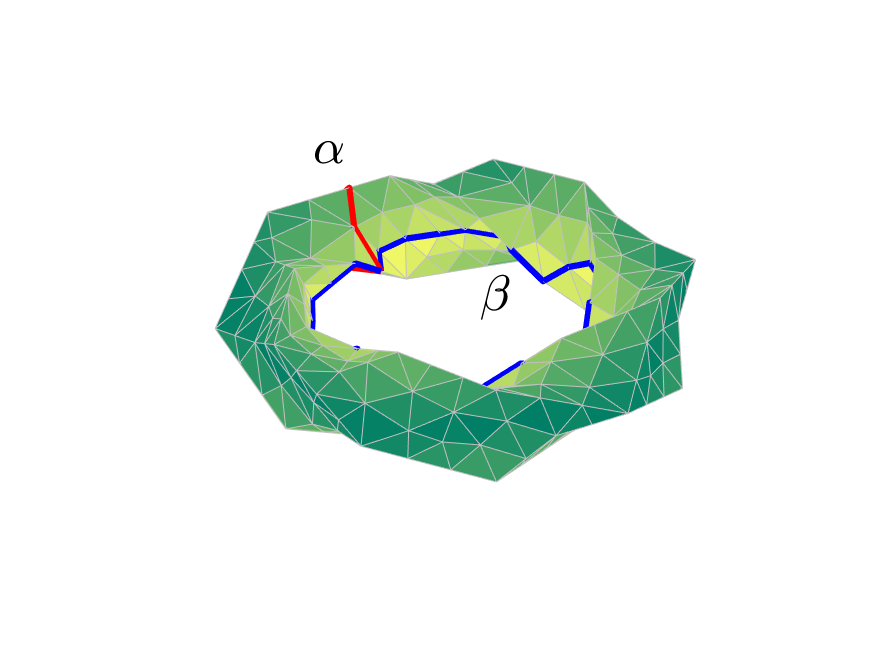}
}
\hspace{0.3cm}
\subfloat[The flattened double periodic domain $\widetilde{\mathcal{T}}$]{\label{subfig:DPplane}
\includegraphics[clip,trim = {1.cm 2cm 0.2cm 1cm},width = 0.4\textwidth]{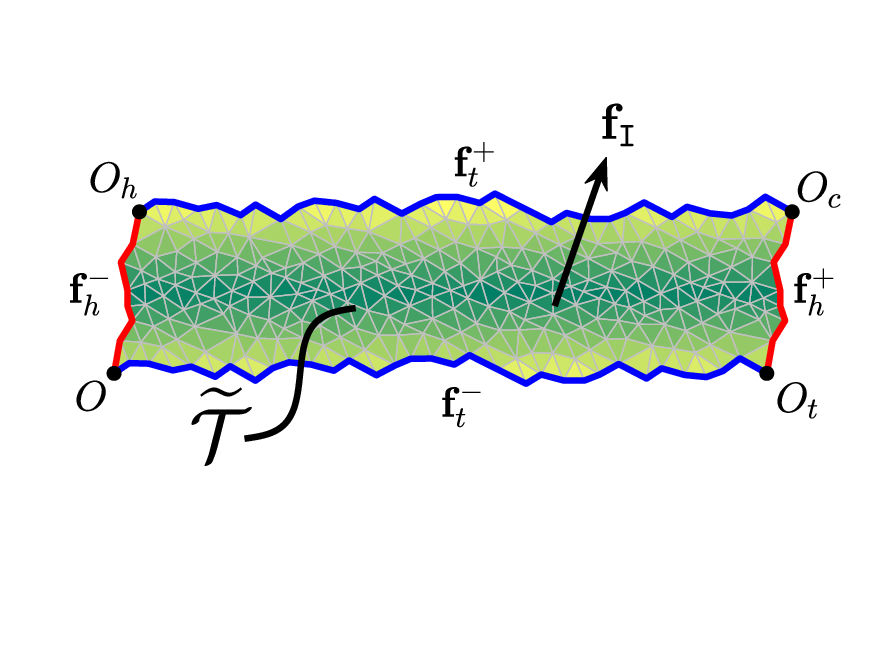}
}
    \caption{An illustrative example of double periodic conformal flattening with four polylined sides $\{\widehat{OO_h}$, $\widehat{O_t O_{c}}$, $\widehat{OO_t}$, $\widehat{O_hO_{c}}\}$. }
    \label{fig:set_classify_T}
\end{figure}

In the discrete version, let $\mathcal{M}$ be a discrete genus-one surface with a handle loop
\begin{subequations} \label{set:loops}
\begin{align}
    \alpha = \{ [a_0,a_1], [a_1,a_2], \cdots, [a_{h-1},a_{h}], [a_h,a_0] \} \subset \mathbb{R}^3, \label{set:handle}
\end{align}
and a tunnel loop
\begin{align}
    \beta = \{ [b_0,b_1], [b_1,b_2], \cdots, [b_{t-1},b_{t}], [b_t,b_0] \} \subset \mathbb{R}^3, \label{set:tunnel}
\end{align}
\end{subequations}
where $a_0 = b_0$, and $\alpha$ and $\beta$ have no self-intersection except at $a_0$, as shown in Figure \ref{subfig:G1sufHT}. We first cut the mesh of $\mathcal{M}$ separately along $\alpha$ and $\beta$ to obtain a new expanded mesh with
\begin{align}
    \widetilde{\mathcal{M}} = \big(\mathcal{M}\setminus (\alpha \cup \beta)\big) \cup \partial \widetilde{\mathcal{M}}, \quad
    \partial \widetilde{\mathcal{M}} = \alpha^+ \cup \alpha^- \cup \beta^+ \cup \beta^-, \label{set:cutM}
\end{align}
where $\alpha = \alpha^+ \# \alpha^-$, $\beta = \beta^+ \# \beta^-$, and `$\#$' denotes the identical gluing operation implemented along the cutting positions.
Let $\widetilde{\mathcal{T}}$ be a double periodic domain with four vertices $\{O = (0,0), O_h = (h^1,h^2), O_t = (t^1,t^2), O_{c} = (h^1+t^1,h^2+t^2)\}$ and four polylined sides $\{ \widehat{OO_h}, \widehat{O_hO_{c}}, \widehat{O_{c}O_t}, \widehat{O_tO} \}$, as in Figure \ref{subfig:DPplane}, to be determined,
in which $\mathbf{f}_\mathtt{I}$ denotes the interior vertices of $\widetilde{\mathcal{T}}$, and $\mathbf{f}_h^\mp$ and $\mathbf{f}_t^\mp$ denote the vertices on the polylined sides $\widehat{OO_h}$, $\widehat{O_t O_{c}}$, $\widehat{OO_t}$ and $\widehat{O_hO_{c}}$, respectively. Let $\mathcal{T}$ be a unit cell of $\Lambda_{h,t}$ that is associated with $\widetilde{\mathcal{T}}$ and satisfies identical conditions
$
     \widehat{O_tO_{c}} = \widehat{OO_h} + t,
     \widehat{O_hO_{c}} = \widehat{OO_t} + h .
$

Let
$\widetilde{\mathcal{M}}$ be the cutting surface of $\mathcal{M}$ as in \eqref{set:loops} along the loops $\alpha$ and $\beta$ as in \eqref{set:loops}. On the basis of the Poincar\'e-Klein-Koebe uniformization theorem, for the genus-one surface $\mathcal{M}$,
we find a conformal parameterization
\begin{subequations} \label{map:MtoT_classfy}
\begin{align}
    \widetilde{f}: \widetilde{\mathcal{M}} \to \widetilde{\mathcal{T}},
\end{align}
with $
    \mathbf{f}:= [\mathbf{f}^1,\mathbf{f}^2] = [\mathbf{f}_\mathtt{I}^\mathrm{T}, (\mathbf{f}_t^-)^\mathrm{T}, (\mathbf{f}_h^-)^\mathrm{T}, O^\mathrm{T}]^\mathrm{T} \in \mathbb{R}^{n\times 2},
$ and
\begin{align}
    \widetilde{\mathbf{f}} = \left[\begin{array}{c}
        \mathbf{f}_\mathtt{I} \\ \mathbf{f}_t^- \\ \mathbf{f}_h^- \\ O \\ \hline \mathbf{f}_t^+ \\ \mathbf{f}_h^+ \\ O_h \\ O_t \\ O_{c}
    \end{array}\right]
    = \left[\begin{array}{cccc|cc}
        I &&&&& \\
        & I &&&& \\
        && I &&& \\
        &&&1&& \\
        \hline
        &I &&& \mathbf{1}& \\
        && I &&& \mathbf{1}\\
        &&&1&1& 0 \\
        &&&1& 0 &1 \\
        &&&1&1&1
    \end{array}\right]
    \left[\begin{array}{c}
        \mathbf{f}_\mathtt{I} \\ \mathbf{f}_t^- \\ \mathbf{f}_h^- \\ O \\ \hline h \\ t
    \end{array}\right]
    := P
    \left[\begin{array}{c}
        \mathbf{f} \\ \hline h \\ t
    \end{array}\right]
    \equiv P\mathbf{g}. \label{eq:ft=Pg}
\end{align}
\end{subequations}

We define the cutting Laplacian matrix for $\widetilde{\mathcal{M}}$ as in \eqref{set:cutM} with
\begin{align}
    \big[ \widetilde{L}_D \big]_{ij} = \begin{cases}
        -\tilde{w}_{ij} \equiv -\frac{1}{2} (\cot \theta_{ij} + \cot \theta_{ji}), & \text{if } [v_{i},v_j] \notin \mathcal{E}(\partial \widetilde{\mathcal{M}}),\\
        -\tilde{w}_{ij} \equiv -\frac{1}{2} \cot \theta_{ij}, & \text{if } [v_{i},v_j] \in \mathcal{E}(\partial \widetilde{\mathcal{M}}),\\
        \sum_{\ell \neq i} \tilde{w}_{i\ell}, & \text{if } i = j, \\
        0, & \text{otherwise},
    \end{cases} \label{eq:LDs}
\end{align}
and consider minimizing the cutting discrete conformal energy functional
\begin{align}
    E_C(\widetilde{f}) = E_D(\widetilde{f}) - A(\widetilde{f})
    = \frac{1}{2} \trace (\widetilde{\mathbf{f}}^\mathrm{~T} \widetilde{L}_D \widetilde{\mathbf{f}})
    - (t^1h^2 - t^2h^1), \label{eq:ECcut}
\end{align}
where $A(\widetilde{f}) = |\widetilde{\mathcal{T}}| = t^1h^2 - t^2h^1$ is the image area of $\widetilde{f}$.
Substituting \eqref{eq:ft=Pg} into \eqref{eq:ECcut}, we get the cutting conformal energy
\begin{align}
    E_C(\mathbf{g})
    &= \frac{1}{2} \trace (\mathbf{g}^\mathrm{T} P^\mathrm{T} \widetilde{L}_D P\mathbf{g})
    - |\widetilde{\mathcal{T}}| \nonumber\\
    &= \frac{1}{2} \trace (\mathbf{g}^\mathrm{T} L \mathbf{g})
    - (t^1h^2 - t^2h^1)
    \equiv E_D(\mathbf{g}) - A(\mathbf{g}), \label{eq:ECcut2}
\end{align}
where $A(\mathbf{g})$ is the image area with respect to $\mathbf{g}$,
\begin{align} \label{eq:PLP}
    L: = P^\mathrm{T} \widetilde{L}_D P = \begin{bmatrix}
        L_D & S \\ S^\mathrm{T} & K
    \end{bmatrix}, \quad
    K = \begin{bmatrix}
        k_{11} & k_{12} \\ k_{12} & k_{22}
    \end{bmatrix}
    \text{symmetric},
\end{align}
$S = [\mathbf{s}_1,\mathbf{s}_2]$ with $\mathbf{1}^\mathrm{T} \mathbf{s}_{j} = 0$, $j = 1,2$, and $L_D$ is the Laplacian matrix corresponding to the vertices on $\mathcal{M}$ as in \eqref{eq:LD}.

We write $\mathbf{g} = [\mathbf{g}^1,\mathbf{g}^2]$, $\mathbf{f} = [\mathbf{f}^1, \mathbf{f}^2]$ and map $f$ from the interacting vertex $a_0$ in \eqref{set:loops} to $O=(0,0)$.
Let the partial derivative of $E_C(\mathbf{g})$ in \eqref{eq:ECcut2} with respect to $\vectx(\mathbf{g})$ be zero. Then
\begin{align}
    \frac{\partial E_C(\mathbf{g})}{\partial \vectx(\mathbf{g})}
    &= \left[ \begin{array}{c}
        L\mathbf{g}^1 \\ L\mathbf{g}^2
    \end{array}
    \right]
    - \left[ \begin{array}{c}
        \mathbf{0} \\ t^2 \\ -h^2 \\ \hline \mathbf{0} \\ -t^1 \\ h^1 \\
    \end{array}
    \right]
     = \left[
    \begin{array}{ccc|ccc}
        L_D & \mathbf{s}_1 & \mathbf{s}_2 &&& \\
        \mathbf{s}_1^\mathrm{T} & k_{11} & k_{12} &&& -1 \\
        \mathbf{s}_2^\mathrm{T} & k_{21} & k_{22} && 1 & \\
        \hline
        &&& L_D & \mathbf{s}_1 & \mathbf{s}_2  \\
        && 1 & \mathbf{s}_1^\mathrm{T} & k_{11} & k_{12}  \\
        & -1 && \mathbf{s}_2^\mathrm{T} & k_{21} & k_{22} \\
    \end{array}
\right]
\left[
    \begin{array}{c}
        \mathbf{f}^1 \\ h^1 \\ t^1 \\
        \hline
        \mathbf{f}^2 \\ h^2 \\ t^2 \\
    \end{array}
\right]
    = 0. \label{eq:pEC=0DPCF}
\end{align}

Since $L_D\mathbf{1} = 0$ and $\mathbf{s}_j^\mathrm{T} \mathbf{1} = 0, j = 1,2$, the homogeneous equation of \eqref{eq:pEC=0DPCF} has nullity $2$. In practice, we set $O = (0,0)$ and $t = (t^1, t^2) = (1,0)$ to be known.
From \eqref{map:MtoT_classfy}, the unknown matrix is $\big[\mathbf{f}_\mathtt{I}^\mathrm{T}, (\mathbf{f}_t^-)^\mathrm{T}, (\mathbf{f}_h^-)^\mathrm{T}, h^\mathrm{T} \big]^\mathrm{T}$. To this end, we let
\begin{subequations}
\begin{align}
    &\hat{L}_0 = L_D(1:end-1,1:end-1),\label{eq:L0} \\
    & \hat{\mathbf{s}}_j = \mathbf{s}_j(1:end-1), j = 1,2,\\
    &\hat{\mathbf{f}} \equiv [\hat{\mathbf{f}}^1, \hat{\mathbf{f}}^2] = \left[ \begin{array}{c} \mathbf{f}_\mathtt{I} \\ \mathbf{f}_t^- \\ \mathbf{f}_h^- \end{array} \right].
\end{align}
\end{subequations}
The derivative of $E_C(\mathbf{g})$ in \eqref{eq:pEC=0DPCF} with respect to
$\left[\begin{array}{c} \hat{\mathbf{f}}^1 \\ h^1 \end{array}\right]$
and
$\left[\begin{array}{c} \hat{\mathbf{f}}^2 \\ h^2 \end{array}\right]$
becomes the following linear system
\begin{align}
    \begin{bmatrix}
        \hat{L}_0 & \hat{\mathbf{s}}_1 \\ \hat{\mathbf{s}}_1^\mathrm{T} & k_{11}
    \end{bmatrix}
    \left[ \begin{array}{c} \hat{\mathbf{f}} \\ h \end{array} \right]
    = \left[ \begin{array}{cc} -\hat{\mathbf{s}}_2 & \mathbf{0} \\ -k_{12} & 1 \end{array} \right]. \label{eq:LS4torus}
\end{align}

Thus, the linear system in \eqref{eq:LS4torus} can be solved to efficiently obtain the desired conformal parameterization,
\begin{subequations} \label{eq:f2ft}
\begin{align}
    \widetilde{\mathbf{f}} = [\mathbf{f}_\mathtt{I}^\mathrm{T}, (\mathbf{f}_t^-)^\mathrm{T}, (\mathbf{f}_h^-)^\mathrm{T}, O^\mathrm{T}, \big|(\mathbf{f}_t^+)^\mathrm{T}, (\mathbf{f}_h^+)^\mathrm{T}, O_h^\mathrm{T}, O_t^\mathrm{T}, O_{c}^\mathrm{T}]^\mathrm{T},
\end{align}
from $\widetilde{\mathcal{M}}$ to $\widetilde{\mathcal{T}}$, for which
\begin{align}
    & \mathbf{f}_t^+ = \mathbf{f}_t^- + (h^1,h^2), \quad \mathbf{f}_h^+ = \mathbf{f}_h^- + (t^1,t^2), \\
    & O_t = (t^1,t^2), \quad O_h = (h^1,h^2), \quad O_{c} = (h^1 + t^1, h^2 + t^2).
\end{align}
\end{subequations}
This approach is summarized in Algorithm \ref{alg:CP-g1}.

\begin{algorithm}[h]
\caption{Conformal parameterization for genus-one surfaces via DPCF}
\begin{algorithmic}[1] \label{alg:CP-g1}
    \REQUIRE A genus-one surface $\mathcal{M}$.
    \ENSURE A conformal flattening map to a double periodic domain $f: \mathcal{M}\to \mathbb{R}^2$.
    \STATE Find a handle loop $\alpha$ and a tunnel loop $\beta$ of $\mathcal{M}$.
    \STATE Cut the surface $\mathcal{M}$ along the paths $\alpha, \beta$ and obtain a simply connected surface $\widetilde{\mathcal{M}}$ with two pairs of identical sides, $\alpha^+, \alpha^-$ and $\beta^+, \beta^-$.
    \STATE Solve linear system \eqref{eq:LS4torus}.
    \STATE Construct the parameterized double periodic domain $\widetilde{\mathcal{T}}$ by the solution $\widetilde{f}$ in \eqref{eq:f2ft}. Then $f: \mathcal{M} \to \mathcal{T}$ is obtained by $\mathbf{f}$ and \eqref{eq:bary}.
\end{algorithmic}
\end{algorithm}

\begin{remark}
Since the Poincar\'e-Klein-Koebe uniformization theorem does not state the relationship between $\omega_1$ and $\omega_2$ in \eqref{group:deckg1}, we do not restrain $h \perp t$ in practical computations but rather automatically compute the optimal relationship between $h$ and $t$ by minimizing the conformal energy instead.
In fact, the relationship between $h$ and $t$ depends on the surface itself and the selected loops.
\end{remark}

\section{Cutting Path Independence} \label{sec:cut}

Now, we aim to prove that the cutting path selection process is independent of the resulting map.
That is, the resulting maps obtained with two pairs of different paths are identical up to cutting and gluing along the edges, global scaling, rotation and translation.
Throughout the following discussion, whenever we mention `identical', it signifies this concept.

\begin{figure}[h]
\centering
\subfloat[The surface $\mathcal{M}$ with cutting paths and $v_*$]{\label{subfig:Mcut}
\includegraphics[clip,trim = {3cm 3cm 2.5cm 0.8cm},width=0.28\textwidth]{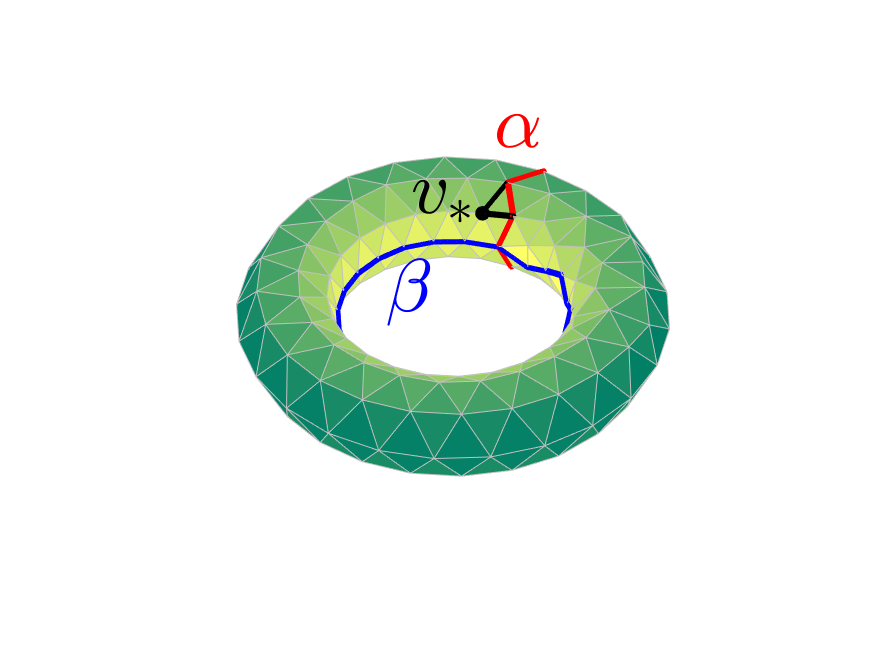}
}
\hspace{0.3cm}
\subfloat[The flattened domain $\widetilde{\mathcal{T}}$ with cutting paths and vertex $\mathbf{f}_*$]{\label{subfig:Tcut}
\includegraphics[clip,trim = {0cm 1.5cm 0cm 0.8cm},width=0.35\textwidth]{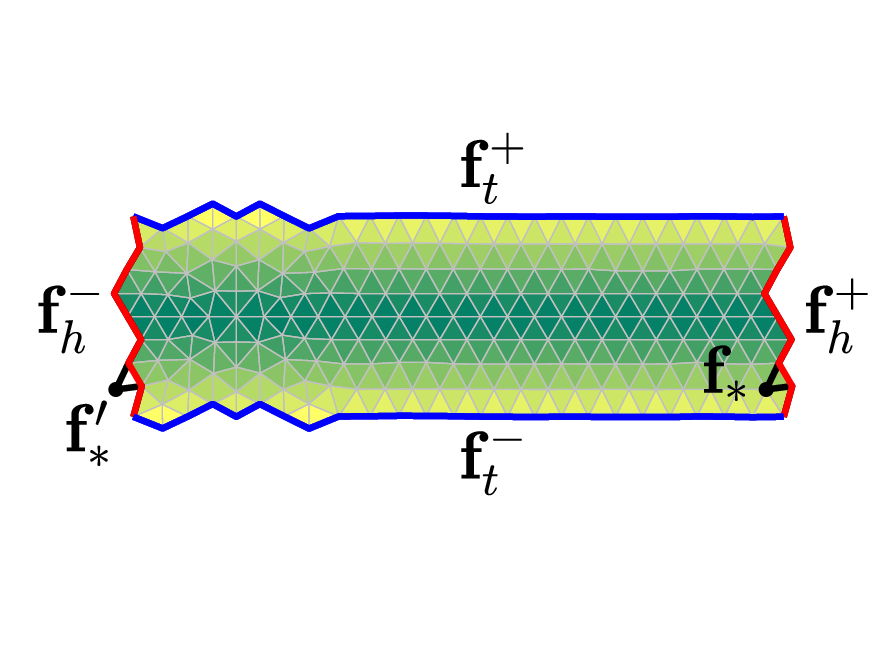}
}
    \caption{The results of DPCF for cutting paths $(\alpha, \beta)$ and $(\alpha', \beta)$, respectively. The left and right ones are cutting paths (red, black and blue lines) on the original surface $\mathcal{M}$ and the flattened domain $\widetilde{\mathcal{T}}$, respectively. }
    \label{fig:cutproof}
\end{figure}

\begin{lemma} \label{lma:cutpath}
    Suppose $\alpha = \{[v_0,v_1], \cdots,  [v_k,v_{k+1}], \cdots, [v_h,v_0]\}$ and $\beta$ be handle and tunnel loop of the given genus-one surface $\mathcal{M}$.
    Let $v_*$ be the adjacent vertex of $\alpha^-$, as in Figure \ref{subfig:Mcut}. Then for the another handle loop $\alpha' = \{[v_0, v_1], \cdots, [v_k, v_*], [v_*, v_{k+1}], \cdots, [v_h, v_0]\}$ with $\beta$.
    Then the two corresponding resulting vertices $f,f'$ and $h,h'$ satisfy
    \begin{align} \label{eq:f'f}
        f'(p) = \begin{cases}
            f(p) - t, & \text{if } p \in [v_k,v_*,v_{k+1}], \\
            f(p), & \text{if } p \notin [v_k,v_*,v_{k+1}],
        \end{cases}
        \quad \text{and} \quad h' = h,
    \end{align}
    with $O' = O = (0,0), t' = t = (1,0)$, as in Figure \ref{subfig:Tcut}.
\end{lemma}

\begin{proof}
Let $\tau = [v_k,v_*,v_{k+1}]$ and separate $\mathbf{f}_* = f(v_*)$ into $\mathbf{f}_*^+$ and $\mathbf{f}_*^-$. Divide $\mathcal{M}$ into $\mathcal{M}\setminus \tau$ and $\tau$
Consider the problem for paths $(\alpha,\beta)$
\begin{subequations}
\begin{align} \label{prob1}
    \min_{f,h} E_C(f,h) \quad \text{s.t.} \ (f,h) \in S(\alpha,\beta):= \left\{(f,h) \left|
    \begin{array}{l}
       \mathbf{f}_i^+ = \mathbf{f}_i^- + h, v_i \in \mathcal{V}(\beta), \\
       \mathbf{f}_i^+ = \mathbf{f}_i^- + t, v_i \in \mathcal{V}(\alpha), \fbox{$\mathbf{f}_*^+ = \mathbf{f}_*^-$}
    \end{array} \right. \right\},
\end{align}
and by $\mathcal{V}(\alpha') = \mathcal{V}(\alpha) \cup \{v_*\}$, the problem for $(\alpha',\beta)$
\begin{align} \label{prob2}
    \min_{f,h} E_C(f,h) \quad \text{s.t.} \ (f,h) \in S(\alpha',\beta):= \left\{(f,h) \left|
    \begin{array}{l}
       \mathbf{f}_i^+ = \mathbf{f}_i^- + h, v_i \in \mathcal{V}(\beta), \\
       \mathbf{f}_i^+ = \mathbf{f}_i^- + t, v_i \in \mathcal{V}(\alpha), \fbox{${\mathbf{f}}_*^+ = {\mathbf{f}}_*^- + t$}
    \end{array} \right. \right\}.
\end{align}
\end{subequations}
For $(\phi,h) \in S(\alpha,\beta)$, define
\begin{align*}
    \mathcal{P}(\phi) = \begin{cases}
            \phi(p) - t, & \text{if } p \in \tau, \\
            \phi(p), & \text{if } p \notin \tau,
        \end{cases}
\end{align*}
and let $\phi' = \mathcal{P}(\phi)$.
Then it can be easily verified that for each pair $(\phi,h) \in S(\alpha,\beta)$, the pair $(\phi',h)$ belongs to $S(\alpha',\beta)$, and vice versa. In other words, the feasible set $S(\alpha,\beta)$ of Problem \eqref{prob1} and the feasible set $S(\alpha',\beta)$ of Problem \eqref{prob2} are homeomorphic under the map $(\mathcal{P},I)$.

On the other hand, from the division $\mathcal{M} = (\mathcal{M}\setminus \tau)\cup\tau$, we have $E_C(f) = E_C(f)\big|_{\mathcal{M}\setminus \tau} + E_C(f)\big|_{\tau}$. Besides, \eqref{eq:f'f} illustrates that
\begin{align} \label{eq:f'ftau}
    \phi'(\tau) = \phi(\tau) - t,\quad \phi'(\mathcal{M}\setminus \tau) = \phi(\mathcal{M}\setminus \tau).
\end{align}
Then we have
\begin{align*}
    E_C(\phi',h) = E_C(\phi',h)\big|_{\mathcal{M}\setminus \tau} + E_C(\phi',h)\big|_{\tau} = E_C(\phi,h)\big|_{\mathcal{M}\setminus \tau} + E_C(\phi - t,h)\big|_{\tau}.
\end{align*}
Since the conformal energy is invariant upon translation, i.e., $E_C(\phi - t,h) = E_C(\phi,h)$, it follows that $E_C(\phi',h) = E_C(\phi,h)$ for each $(\phi,h) \in S(\alpha,\beta)$ and $(\mathcal{P}(\phi),h)\in S(\alpha',\beta)$.
Finally, $(f,h)$ is the minimizer of Problem \eqref{prob1} if and only if $(f',h) = (\mathcal{P}(f),h)$ is the minimizer of Problem \eqref{prob2}, which completes the proof.

\end{proof}

This lemma also holds when the tunnel loop path $\beta$ changes.
Intuitively, this lemma illustrates that if we modify the path so that a triangle changes from one side of the path to the other, the resulting flattened domain can be obtained by cutting the corresponding triangle and gluing it to the another side.
Then, we can prove the cutting path independence of our proposed DPCF algorithm.

\begin{theorem} \label{thm:cutglue}
    Given a genus-one triangular mesh $\mathcal{M}$ and two pairs of cutting paths $(\alpha,\beta)$ and $(\alpha',\beta')$, the corresponding two maps $f$ and $f'$ obtained by DPCF in Algorithm \ref{alg:CP-g1} are identical up to cutting and gluing along the edges, global scaling, rotation and translation.
\end{theorem}

\begin{proof}
    Upon repeatedly using Lemma \ref{lma:cutpath} to move loops, the proof is completed.
\end{proof}

This theorem explains that the maps resulting from DPCF are seamless, which is an advantage of our algorithm: the selected cutting path does not affect the result. Thus, we do not have to incur an additional cost by seeking `optimal' paths to reduce the conformal distortion.
Following this theorem, the relationship between cutting paths $\alpha,\beta$ and $h,t$ can be immediately obtained by cutting and gluing the flattened triangulation.
\begin{theorem}
    Let $(\omega_1,\omega_2)$ and $(\omega'_1,\omega'_2)$ be two homology bases corresponding to two pairs of cutting paths $(\alpha,\beta)$ and $(\alpha',\beta')$, respectively, and satisfy
    \begin{align*}
        \begin{bmatrix}
            \omega'_1 \\ \omega'_2
        \end{bmatrix}
        =
        \begin{bmatrix}
            a_{11} & a_{12} \\ a_{21} & a_{22}
        \end{bmatrix}
        \begin{bmatrix}
            \omega_1 \\ \omega_2
        \end{bmatrix}, \quad a_{ij} \in \mathbb{Z}, i,j = 1,2.
    \end{align*}
    Then, the resulting translation vectors $(h,t)$ and $(h',t')$ produced by DPCF in Algorithm \ref{alg:CP-g1} satisfy
    \begin{align*}
        \begin{bmatrix}
            h' \\ t'
        \end{bmatrix}
        =
        \begin{bmatrix}
            a_{11} & a_{12} \\ a_{21} & a_{22}
        \end{bmatrix}
        \begin{bmatrix}
            h \\ t
        \end{bmatrix} G, \quad G^\mathrm{T} G = \lambda I_2, \lambda>0 \text{ is a constant.}
    \end{align*}
\end{theorem}

\begin{figure}
    \centering
    \begin{tabular}{m{0.3\textwidth}<{\centering}m{0.3\textwidth}<{\centering}m{0.3\textwidth}<{\centering}}
        \includegraphics[clip,trim = {1cm 3cm 1cm 2.3cm},width=\linewidth]{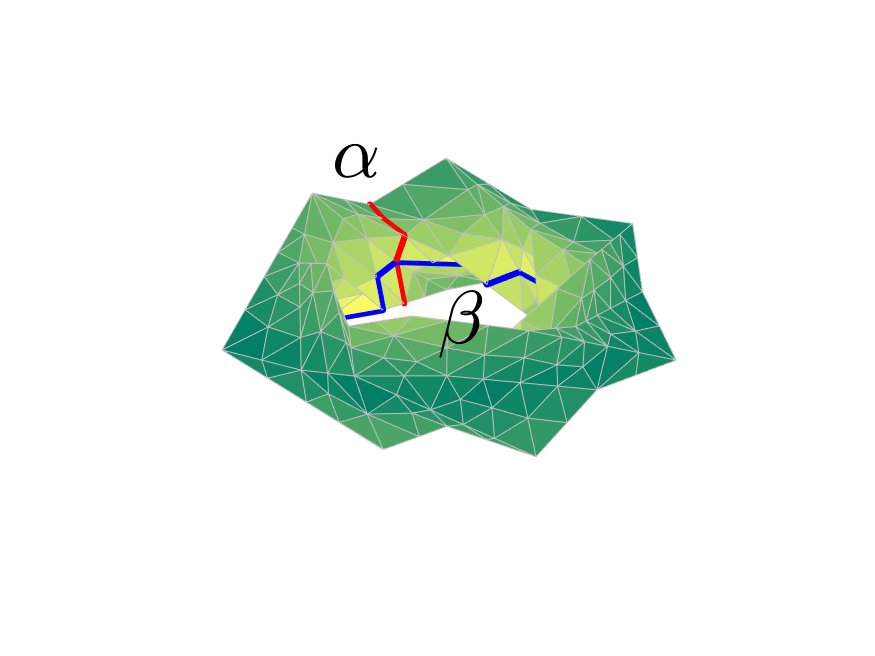} &
        \includegraphics[clip,trim = {1cm 3cm 1cm 2.3cm},width=\linewidth]{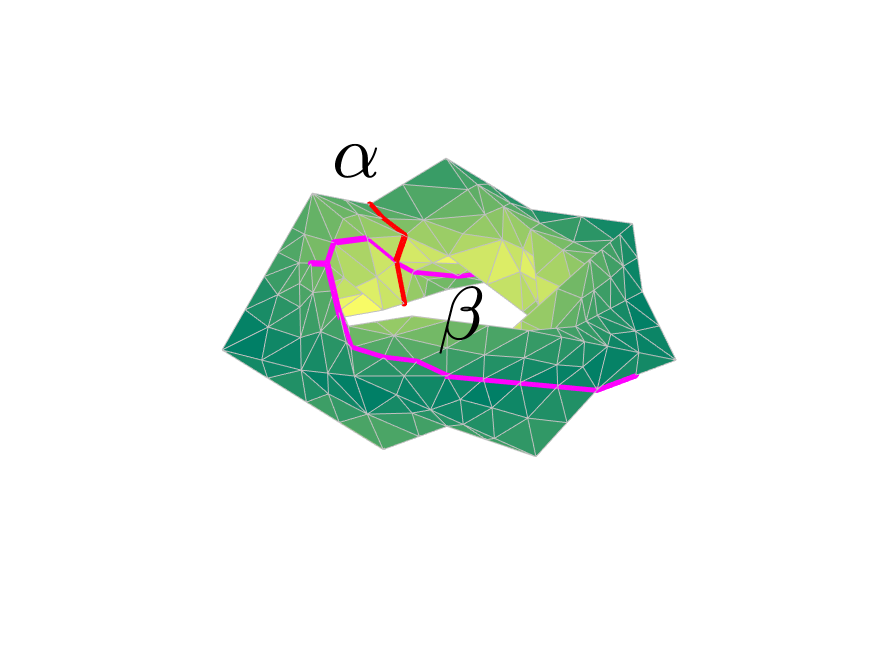} &
        \includegraphics[clip,trim = {1cm 3cm 1cm 2.3cm},width=\linewidth]{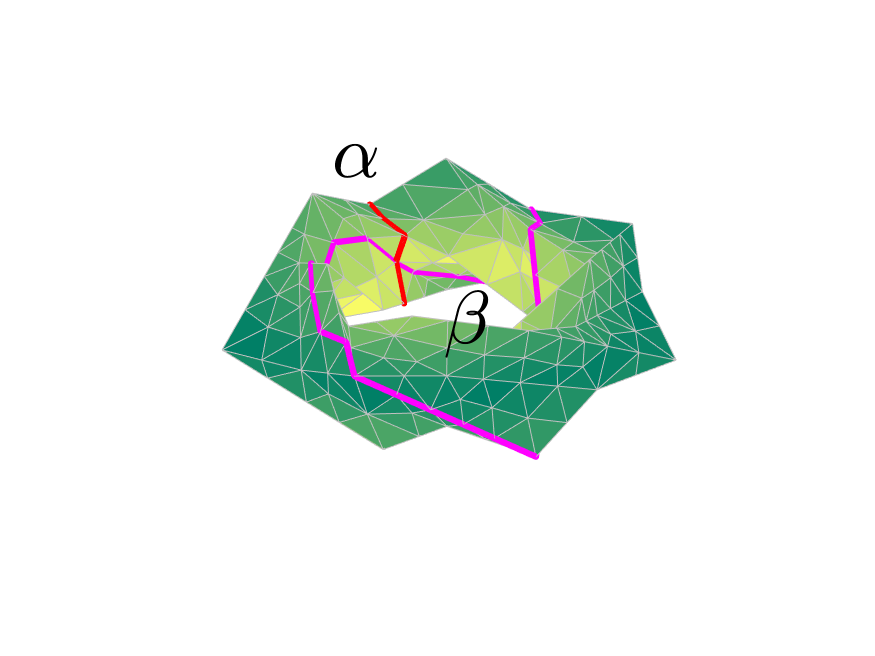} \\
        \includegraphics[clip,trim = {0.5cm 1.5cm 0.5cm 0.8cm},width=\linewidth]{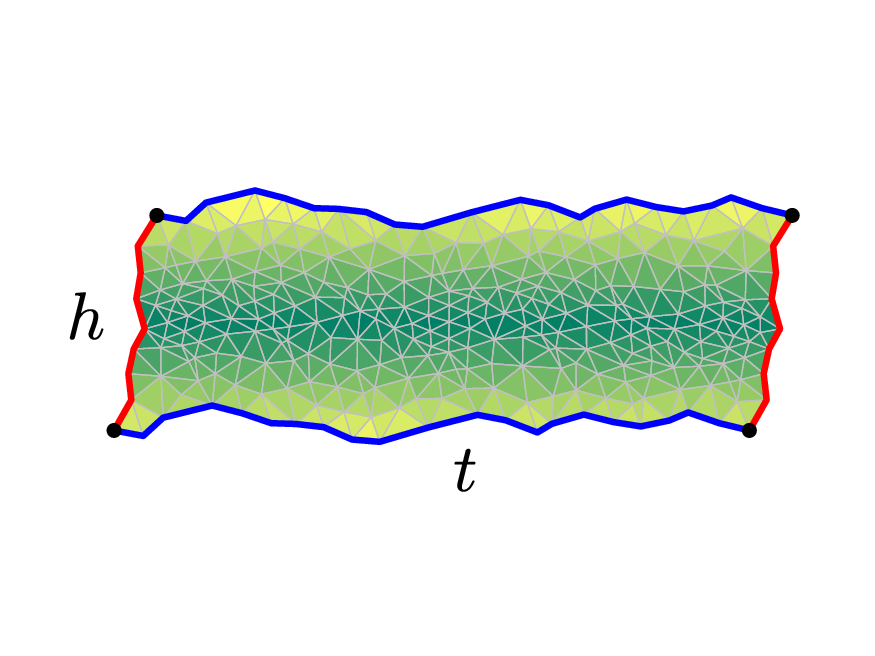} &
        \includegraphics[clip,trim = {0.5cm 1.5cm 0.8cm 0.8cm},width=\linewidth]{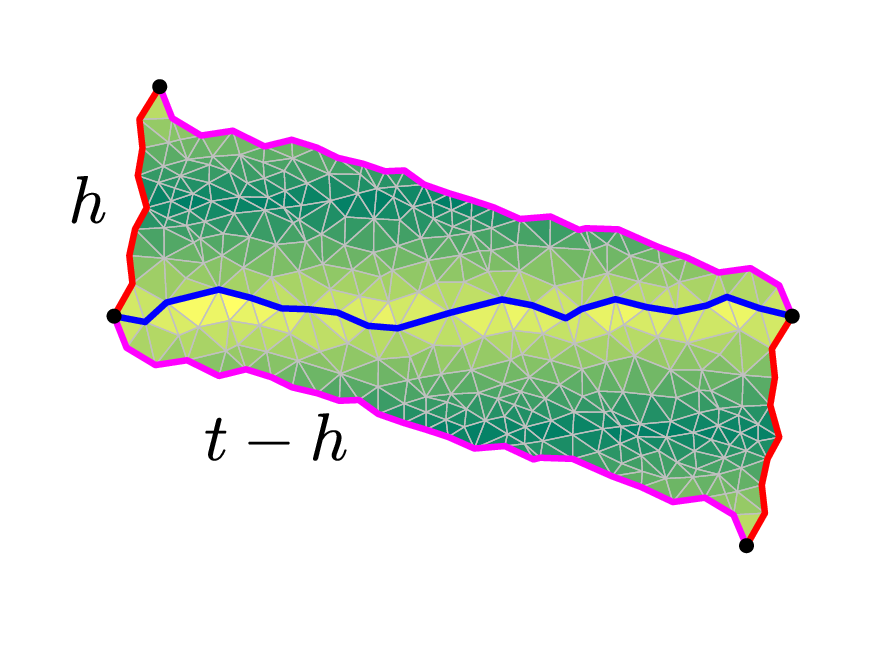} &
        \includegraphics[clip,trim = {1.3cm 0.8cm 0.8cm 0.5cm},width=\linewidth]{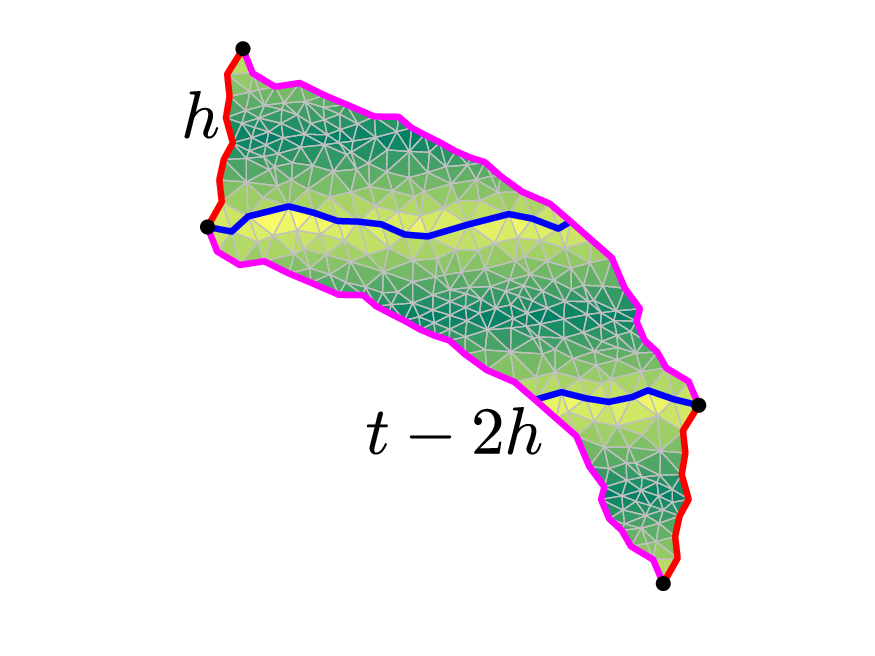} \\
        \small (a) Cutting paths belonging to homology basis $\omega_1,\omega_2$. &
        \small (b) Cutting paths belonging to homology basis $\omega_1,\omega_2-\omega_1$. &
        \small (c) Cutting paths belonging to homology basis $\omega_1,\omega_2-2\omega_1$.
    \end{tabular}
    \caption{A genus-one surface and flattened domains with cutting paths belonging to homology basis $\omega_1,\omega_2$, $\omega_1,\omega_2-\omega_1$ and $\omega_1,\omega_2-2\omega_1$, respectively. The red and blue paths are the loops belonging to homology basis $\omega_1,\omega_2$, and the magenta paths are the tunnel loops on belonging to $\omega_2 - \omega_1$ and $\omega_2 - 2\omega_1$, respectively.  }
    \label{fig:cut&glue}
\end{figure}

\begin{proof}
    For different cutting paths, Theorem \ref{thm:cutglue} states that the flattened domains are identical. Hence, we can cut the domain along the cutting paths and glue them along the cutting seams accordingly, as in Figure \ref{fig:cut&glue}.
    Since we set $O = (0,0)$ and $t = (1,0)$ in Algorithm \ref{alg:CP-g1}, there may exist a rotation and scaling transformations between the translation vectors $(h,t)$ and $(h',t')$ derived from different cutting paths.
    Thus, the proof is completed.
\end{proof}

This theorem illustrates that the translation vectors $(h,t)$ only relate to the homology bases to which the cutting paths belong.
Hence, we can control the relationship between $h$ and $t$ by cutting and gluing the flattened triangulation. The corresponding cutting paths are obtained simultaneously from the boundary of the new triangulation.

\section{Bijectivity Guarantee} \label{sec:bijective}

The bijectivity of the conformal map is not always guaranteed.
Gortler et al. \cite{SGCG06} proved the bijectivity for the 1-form methods.
Floater \cite{Floa02} proved that the harmonic map is bijective under the assumption that each edge weight $w_{ij}$ is positive and the image domain is convex.
However, this result cannot be directly applied for DPCF since (a) the convexity is not satisfied under periodic boundary constraints, and (b) the vectors $(h,t)$ should be obtained simultaneously while compute the conformal map $\mathbf{f}$, such that the domain is not prescribed.
In this case, we generalize the result of \cite{Floa02} and prove that the bijectivity of DPCF is still guaranteed if the weight $w_{ij}$ is positive for each edge.

Under the positive edge weight assumption, we have the following theorem concerning the flattened domain.
\begin{theorem} \label{thm:h20}
For the genus-one triangular mesh $\mathcal{M}$, if the edge weight $w_{ij}$ in \eqref{eq:LD} is positive, then the solution of \eqref{eq:LS4torus} satisfies $A(\widetilde{f}) = h^2>0$.
\end{theorem}

\begin{proof}
    Given $O = (0,0)$ and $t = (1,0)$, the second system of \eqref{eq:LS4torus} is
    \begin{align*}
        \begin{bmatrix}
        \hat{L}_0 & \hat{\mathbf{s}}_1 \\ \hat{\mathbf{s}}_1^\mathrm{T} & k_{11}
    \end{bmatrix}
    \left[ \begin{array}{c} \hat{\mathbf{f}}^2 \\ h^2 \end{array} \right]
    = \left[ \begin{array}{c} \mathbf{0} \\ -1 \end{array} \right].
    \end{align*}
    If $h^2 = 0$, we have $\hat{L}_0\hat{\mathbf{f}}^2 = 0$. Since $\hat{L}_0$ is an M-matrix, then  $\hat{\mathbf{f}}^2 = 0$, contradicting to $\mathbf{s}_1^\mathrm{T} \hat{\mathbf{f}}^2 + k_{11} h^2 = -1$. Hence, $h^2 \neq 0$.

    Assuming that $\mathbf{g}$ with $h^2<0$ is the solution of \eqref{eq:LS4torus}, we have $A(\mathbf{g}) = t^1h^2 - t^2h^1 = h^2 < 0$ for $t = (1,0)$.
    Define $\mathbf{g}' = [\mathbf{g}^1,-\mathbf{g}^2] = \mathbf{g}\begin{bmatrix} 1 & \\ & -1 \end{bmatrix}$. We can verify that $E_D(\mathbf{g}') = E_D(\mathbf{g})$ and $A(\mathbf{g}') = -A(\mathbf{g}) > A(\mathbf{g})$ in \eqref{eq:ECcut2}. Hence, we have $E_C(\mathbf{g}') < E_C(\mathbf{g})$, which contradicts the statement that $E_C(\mathbf{g})$ reaches the minimum. Hence, we obtain $h^2 > 0$.
\end{proof}

In this theorem, the positive edge weight condition is not necessary. As long as the matrix $\hat{L}_0$ is invertible, we can prove $h^2>0$. Theorem \ref{thm:h20} illustrates that the fundamental domain $\widetilde{\mathcal{T}}$ by DPCF in Algorithm \ref{alg:CP-g1} cannot degenerate to a polyline and is a 2D domain with positive directed area $A(\widetilde{f}) = h^2 > 0$.

\begin{figure}
    \centering
    \subfloat[]{\label{subfig:f1234}
    \includegraphics[clip,trim = {1.cm 2.8cm 0.cm 2.cm},width=0.45\textwidth]{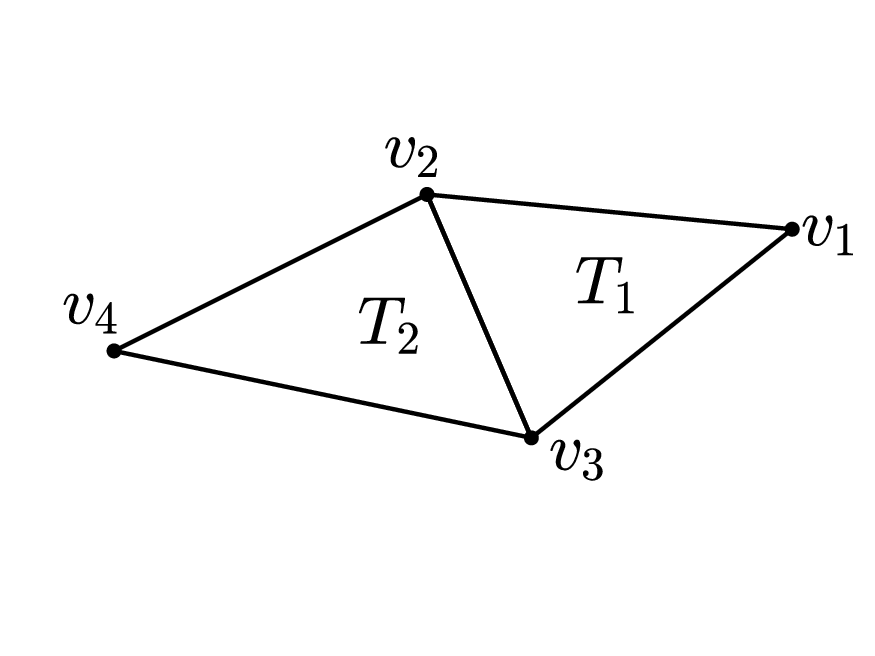}
    }
    \hspace{0.2cm}
    \subfloat[]{\label{subfig:f1234_invert}
    \includegraphics[clip,trim = {2cm 3.2cm 0.5cm 2.5cm},width=0.45\textwidth]{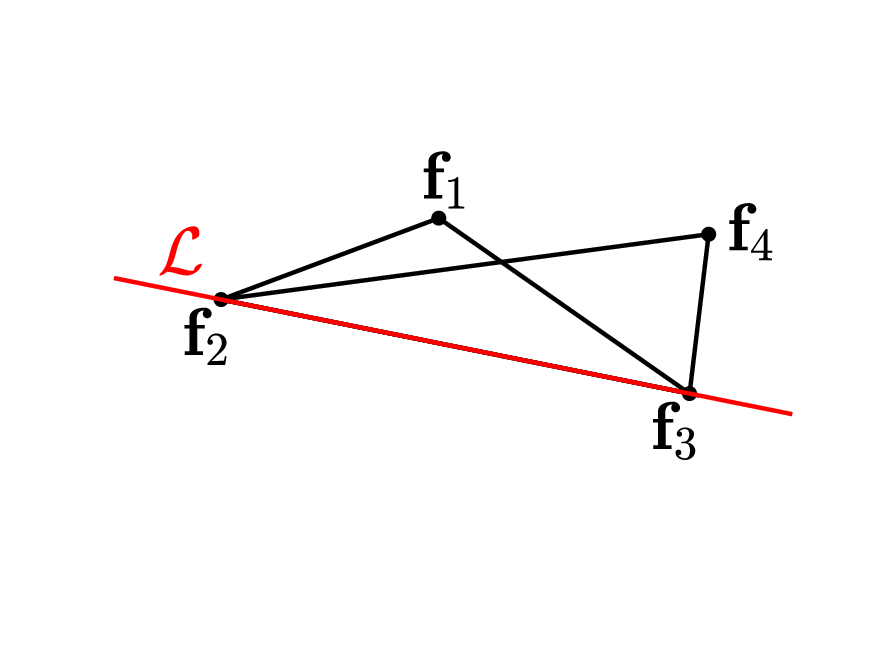}
    }
    \caption{(a) $T_1\cup T_2$ is a quadrilateral in $\widetilde{\mathcal{M}}$. (b) $\mathbf{f}_1$ and $\mathbf{f}_4$ lie on the same side of $\mathcal{L}$ passing through $\mathbf{f}_2$ and $\mathbf{f}_3$.}
    \label{fig:f1234}
\end{figure}

If the edge weight $w_{ij}$ in the discrete harmonic equation is positive, then the resulting map $\widetilde{f}$ is a convex combination map \cite{Floa02} satisfying
\begin{align}
    \widetilde{\mathbf{f}}_i = \sum_{j \in \mathcal{N}(i)} \lambda_{ij} \widetilde{\mathbf{f}}_j, \quad \lambda_{ij} = \frac{w_{ij}}{\sum_{j \in \mathcal{N}(i)} w_{ij}} > 0, \quad \sum_{j \in \mathcal{N}(i)} \lambda_{ij} = 1,  \label{eq:fConvexCombination}
\end{align}
for each interior vertex $v_i$.
For the convex combination function, we have the associated discrete maximum principle \cite{Floa02}.
\begin{lemma}[Discrete Maximum Principle \cite{Floa02}] \label{lma:DMP}
    Let $f$ be a convex combination function over a triangulation $\mathcal{M}$. For any interior vertex $v_0$ of $\mathcal{M}$, let $V_0$ denote the set of all boundary vertices which can be connected to $v_0$ by an interior path. If $f(v_0)\geq f(v)$ for all $v\in V_0$, then $f(v) = f(v_0)$ for all $v\in V_0$.
\end{lemma}
In the following lemma, we want to show that $\widetilde{f}$ is locally bijective based on the discrete maximum principle \cite{Floa02}.
Thus, in light of Theorems 4.1 and 4.2 from \cite{Floa02}, we can prove that the map $\widetilde{f}$ is bijective over $\widetilde{\mathcal{M}}$.
\begin{lemma} \label{lma:bijective}
Suppose that the genus-one surface $\mathcal{M}$ has positive edge weight $w_{ij}$ in \eqref{eq:LD}. Let $T_1 \cup T_2$ be a quadrilateral in $\widetilde{\mathcal{T}}$, as in Figure \ref{subfig:f1234}. If $\widetilde{f}\big|_{T_1}$ is bijective, then $\widetilde{f}\big|_{T_1\cup T_2}$ is also bijective, i.e., $\mathbf{f}_1$ and $\mathbf{f}_4$ cannot be located on the same side of the line $\mathcal{L}$ passing through $\mathbf{f}_2$ and $\mathbf{f}_3$ as in Figure \ref{subfig:f1234_invert}.
\end{lemma}

\begin{figure}
    \centering
    \includegraphics[width=0.75\textwidth]{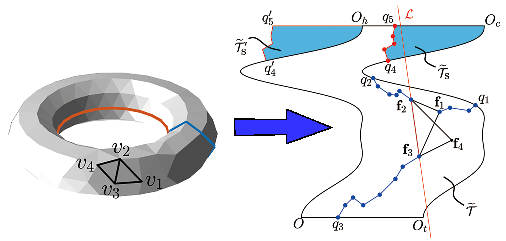}
    \caption{The genus-one surface $\mathcal{M}$ and its fundamental domain $\widetilde{\mathcal{T}}(O,O_t,O_{c},O_h)$ under the convex combination map $\widetilde{f}$. }
    \label{fig:BijectiveProof}
\end{figure}

\begin{proof}
Since $\mathcal{M}$ has positive edge weight $w_{ij}$, the flattened domain $\widetilde{\mathcal{T}}$ does not degenerate to a polyline according to Theorem \ref{thm:h20}.
    Let $\ell$ with
    \begin{align*}
        \ell(\mathbf{f}) = a\mathbf{f}^1 + b\mathbf{f}^2 +c, \quad \mathbf{f} = (\mathbf{f}^1,\mathbf{f}^2)
    \end{align*}
    be the straight line $\mathcal{L}$ passing through $\mathbf{f}_2$ and $\mathbf{f}_3$ with $\ell(\mathbf{f}_2) = \ell(\mathbf{f}_3) = 0$.
    Suppose that $\mathbf{f}_1$ and $\mathbf{f}_4$ lie on the same side of $\mathcal{L}$, i.e., $\ell(\mathbf{f}_1)> 0$ and $\ell(\mathbf{f}_4)\geq 0$ without loss of generality.
    In fact, the fundamental domain $\widetilde{\mathcal{T}}$ in Figure \ref{fig:BijectiveProof} and $T_1 \cup T_2$ with images $\mathbf{f}_1,\mathbf{f}_2,\mathbf{f}_3,\mathbf{f}_4$ in Figure \ref{subfig:f1234_invert} can be used to claim the general case.

    From the case presented in Figure \ref{fig:BijectiveProof}, $\ell < 0$ on the left side of $\mathcal{L}$, and $\ell>0$ on the right side of $\mathcal{L}$.
    We define $H = \ell\circ \widetilde{f}: \widetilde{\mathcal{M}} \to \mathbb{R}$.
    Since $\widetilde{f}$ is a convex combination map, $H$ is a convex combination function with $H(v_2) = H(v_3) = 0$, $H(v_4) \geq 0$ and $H(v_1) > 0$.
    Then, we can find a rising path from $\mathbf{f}_1$ and two falling paths from $\mathbf{f}_2$ and $\mathbf{f}_3$, which hit the boundary of $\widetilde{\mathcal{T}}$ at $q_1:=\widetilde{f}(p_1),q_2:=\widetilde{f}(p_2),q_3:=\widetilde{f}(p_3)$, where $\{p_i\}_{i=1}^3 \subseteq \mathcal{V}(\widetilde{\mathcal{M}})$, along the edges of $\mathcal{E}(\widetilde{\mathcal{T}})$, respectively, as in Figure \ref{fig:BijectiveProof}.
    According to the discrete maximum principle in Lemma \ref{lma:DMP}, paths $\widehat{v_2p_2} = \widetilde{f}^{-1}(\widehat{\mathbf{f}_2q_2})$ and $\widehat{v_3p_3} = \widetilde{f}^{-1}(\widehat{\mathbf{f}_3q_3})$ cannot intersect at a common vertex \cite{Floa02}.

    On the polygon domain $\mathcal{R}$ with
    \[
    \partial \mathcal{R} = \widehat{\mathbf{f}_2q_2} \cup \widehat{q_2O_{c}} \cup \widehat{O_{c}O_h} \cup \widehat{O_hO} \cup \widehat{Oq_3} \cup \widehat{q_3\mathbf{f}_3} \cup \widehat{\mathbf{f}_3\mathbf{f}_2},
    \]
    we can see that $\ell(f) < 0$ except a domain near $O_{c}$.
    For this domain, we can find a subtriangulation domain $\widetilde{\mathcal{T}}_{\text{s}} \subset \widetilde{\mathcal{T}}$ with boundary $\widehat{q_4O_{c}} \cup \widehat{O_{c}q_5} \cup \widehat{q_5q_4}$, which satisfies $\ell(\mathbf{f}_i)<0, \mathbf{f}_i \in \widehat{q_5q_4}$.
    We now perform cutting $\widetilde{\mathcal{T}}_{\text{s}}$ and gluing to $\widetilde{\mathcal{T}}_{\text{s}}'$ with $\widetilde{\mathcal{T}}_{\text{s}}' = \widetilde{\mathcal{T}}_{\text{s}} - t$, whose boundary is $\widehat{q_4'O_h} \cup \widehat{O_hq_5'} \cup \widehat{q_5'q_4'}$.
    Finally, we find a new polygon domain $\mathcal{R}'$ with
    \[
    \partial \mathcal{R}' = \widehat{\mathbf{f}_2q_2} \cup \widehat{q_2q_4} \cup \widehat{q_4q_5} \cup \widehat{q_5O_h} \cup \widehat{O_hq_5'} \cup \widehat{q_5'q_4'} \cup \widehat{q_4'O} \cup \widehat{Oq_3} \cup \widehat{q_3\mathbf{f}_3} \cup \widehat{\mathbf{f}_3\mathbf{f}_2},
    \]
    for which we have $\ell(\mathbf{f}_i) < 0, \mathbf{f}_i \in \partial \mathcal{R}'$. If $\ell(\mathbf{f}_i) < 0$ is not satisfied, we can repeatedly cut and glue the subtriangulations until it is satisfied.
    Correspondingly, for the new fundamental domain $\widetilde{\mathcal{T}}':=\widetilde{\mathcal{T}} - \widetilde{\mathcal{T}}_{\text{s}} + \widetilde{\mathcal{T}}_{\text{s}}'$, we can obtain new loops, cutting surface $\widetilde{\mathcal{M}}'$ and a flattening map $\widetilde{f}': \widetilde{\mathcal{M}}' \to \widetilde{\mathcal{T}}'$ satisfying \eqref{eq:fConvexCombination} for the interior vertices according to Theorem \ref{thm:cutglue}.
    Therefore, $\widetilde{f}'$ is a convex combination map, so it follows that $H' = \ell\circ \widetilde{f}'$ is also a convex combination function.
    Clearly, we have $H'(v_i)<0, v_i \in \widetilde{f}'^{-1} (\partial \mathcal{R}')$.
    By the discrete maximum principle in Lemma \ref{lma:DMP}, we have $H'(v_i)<0, v_i \in \widetilde{f}'^{-1} (\mathcal{R}')$. Since $v_4$ is enclosed by $\widetilde{f}'^{-1} (\partial \mathcal{R}')$, we also have $H'(v_4) = \ell(\mathbf{f}_4) < 0$, which yields a contradiction.
\end{proof}

Now, the bijectivity can be guaranteed.
\begin{theorem} \label{thm:delG1}
    For the genus-one triangular mesh $\mathcal{M}$, if the edge weight $w_{ij}$ in \eqref{eq:LD} is positive, the map $f$ resulting from DPCF in Algorithm \ref{alg:CP-g1} is bijective.
\end{theorem}

\begin{proof}
    Based on Theorem 4.2, Lemma 5.3 and 5.4 in \cite{Floa02} and Lemma \ref{lma:bijective}, the proof is completed.
\end{proof}

\section{Numerical Experiments and Applications} \label{sec:NE}

In this section, we describe the numerical performance and phenomena of our proposed algorithm from various perspectives. All the experimental routines are executed in MATLAB R2024a on a personal computer with a 2.50 GHz CPU and 64 GB RAM. Most of the mesh models are obtained from Thingi10K \cite{Thingi10K} and Common 3D Test Models (\url{https://github.com/alecjacobson/common-3d-test-models}), some of which are remeshed properly. Some special meshes are manually generated by MATLAB built-in function \emph{generateMesh} and proper deformations.
The handle and tunnel loops are computed by ReebHanTun \cite{TDFF13}.
We use two terms to measure the conformal accuracy.
\begin{itemize}
    \item The absolute angle error $\delta$,
    \begin{align}
        \delta(\theta_{ij}) = |\theta_{ij} - \theta_{ij}(f)|, \label{err:aae}
    \end{align}
    where $\theta_{ij}(f)$ means the angle on $\widetilde{\mathcal{T}}$ corresponding to the angle $\theta_{ij}$ on $\mathcal{M}$.
    $\delta$ measures the absolute error of each angle upon the map $f$.

    \item The Beltrami coefficient $\mu$ \cite{PTKC15},
    \begin{align}
        \mu = \frac{\partial f}{\partial \bar{z}} \left/ \frac{\partial f}{\partial z}\right., \label{err:Beltrami}
    \end{align}
    where $\frac{\partial}{\partial \bar{z}} = \frac{1}{2}(\frac{\partial }{\partial x} + \mathrm{i}\frac{\partial }{\partial y})$, $\frac{\partial}{\partial z} = \frac{1}{2}(\frac{\partial }{\partial x} - \mathrm{i}\frac{\partial }{\partial y})$, $z = x + \mathrm{i}y$, and $\mathrm{i} := \sqrt{-1}$ is the imaginary unit.
    The modulus $|\mu|$ measures the proximity of $f$ to the conformality. The map $f$ is conformal if and only if $|\mu| = 0$.
\end{itemize}

\subsection{Double Periodic Conformal Flattening}
\begin{figure}[htp]
    \centering
    \resizebox{\textwidth}{!}{
    \begin{tabular}{cccccc@{}c}
        \includegraphics[clip,trim = {4cm 1cm 3.5cm 0.8cm},width = 0.15\textwidth]{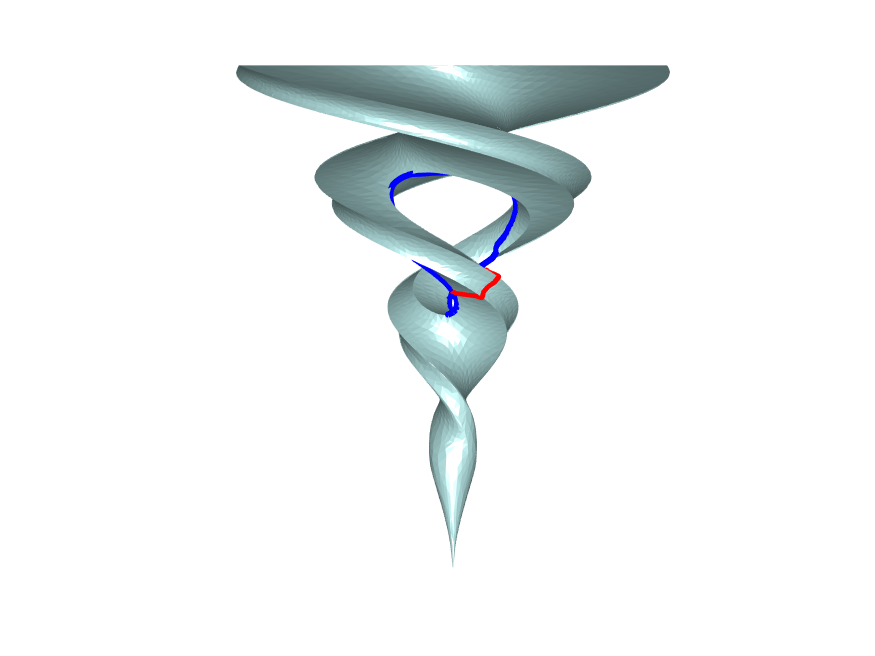} &
        \includegraphics[clip,trim = {4cm 1cm 3.8cm 0.8cm},width = 0.15\textwidth]{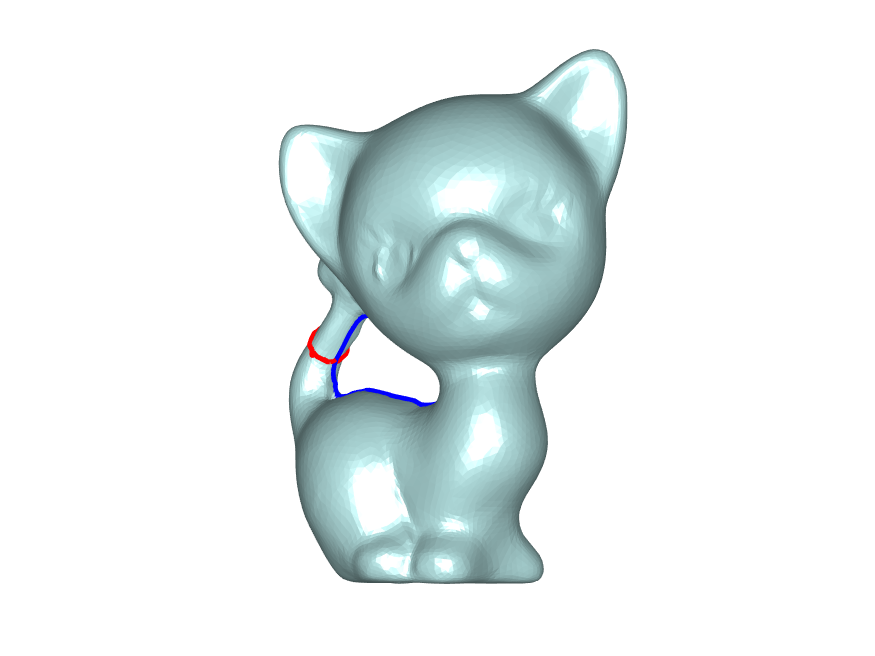} &
        \includegraphics[clip,trim = {4cm 1cm 3.8cm 0.8cm},width = 0.15\textwidth]{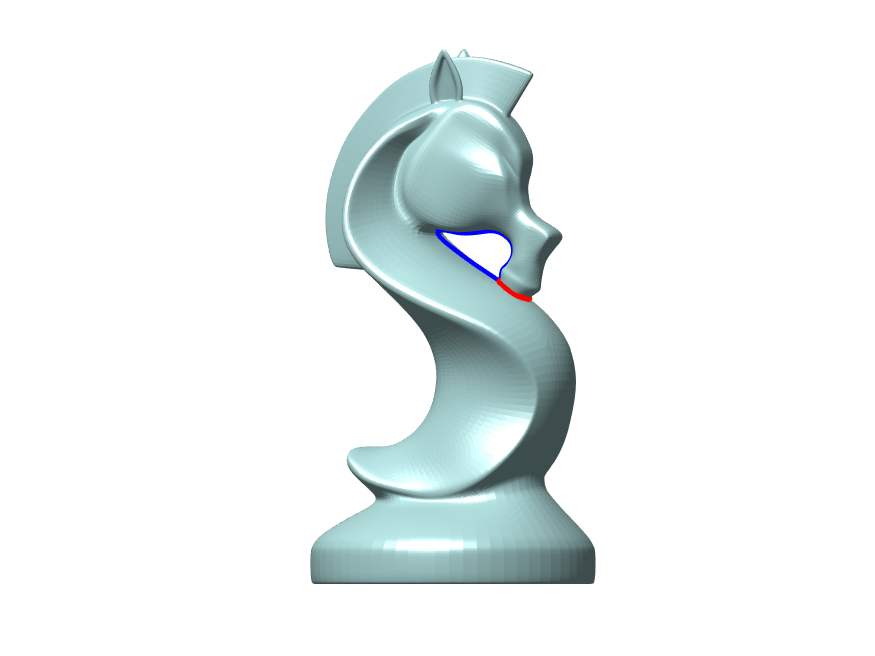} &
        \includegraphics[clip,trim = {4cm 2cm 3.5cm 0.8cm},width = 0.15\textwidth]{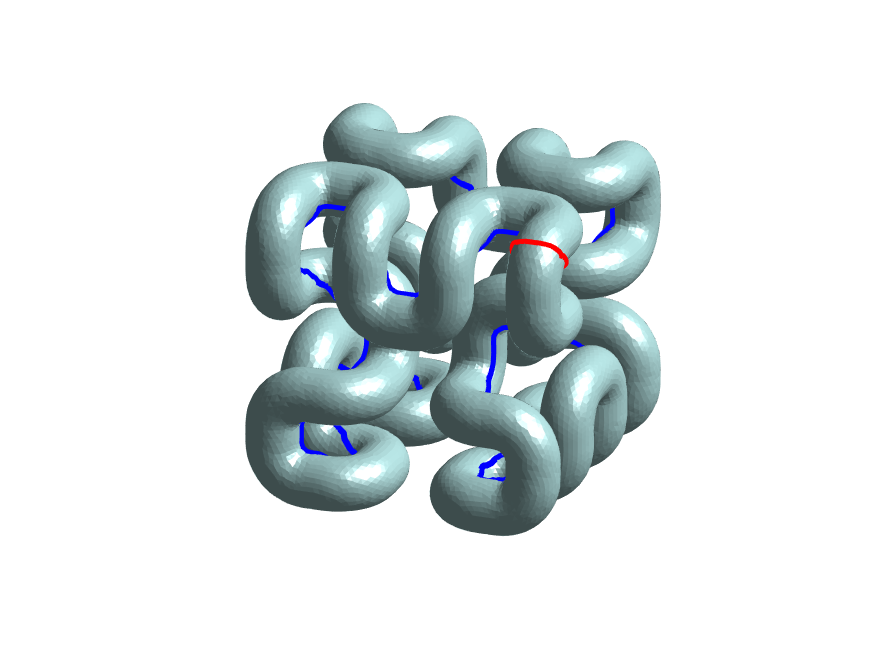} &
        \includegraphics[clip,trim = {3.cm 1cm 2.8cm 1cm},width = 0.15\textwidth]{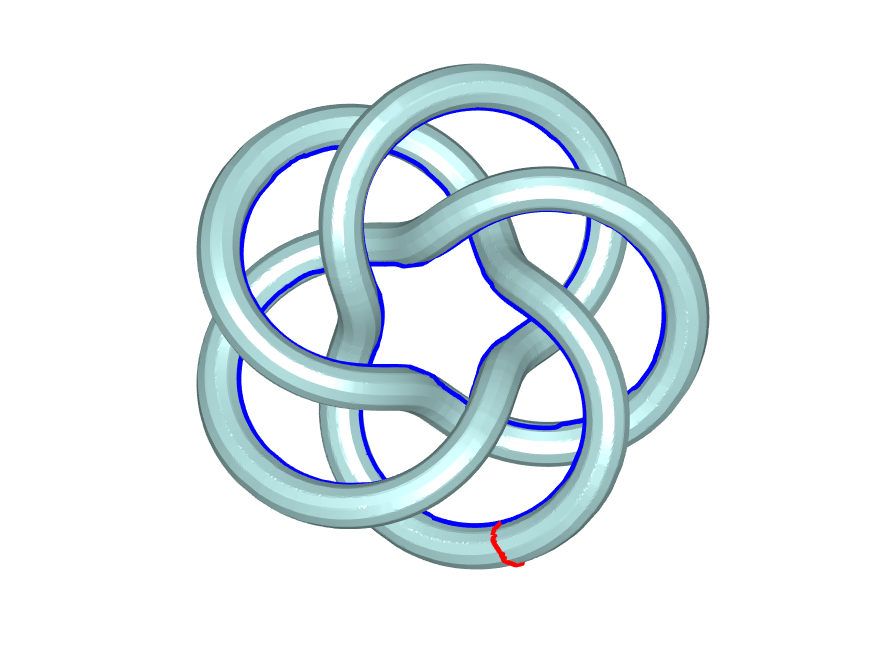} &
        \includegraphics[clip,trim = {5cm 2cm 3.8cm 2cm},width = 0.15\textwidth]{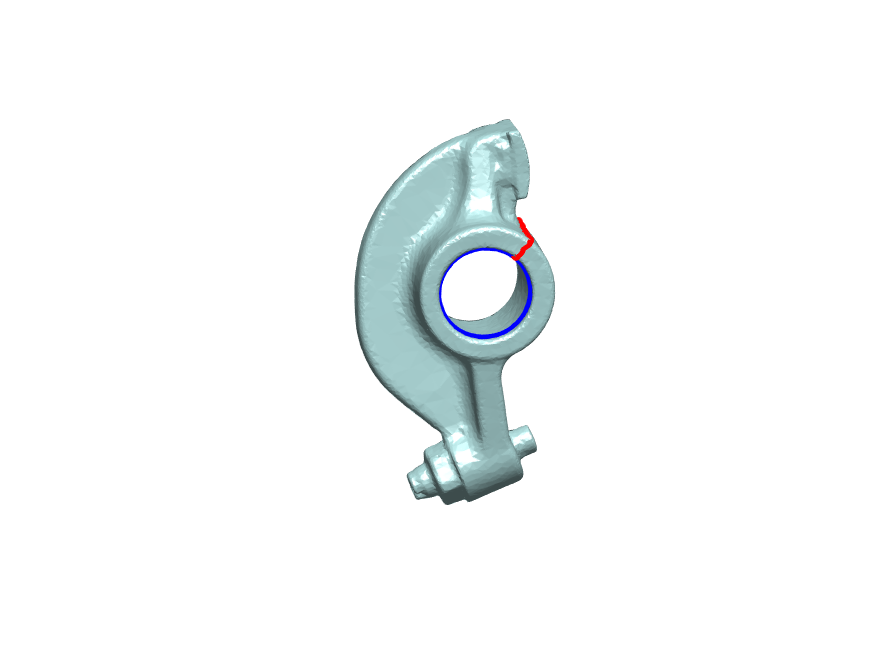} & \\
        \includegraphics[clip,trim = {4cm 1cm 3.5cm 0.8cm},width = 0.15\textwidth]{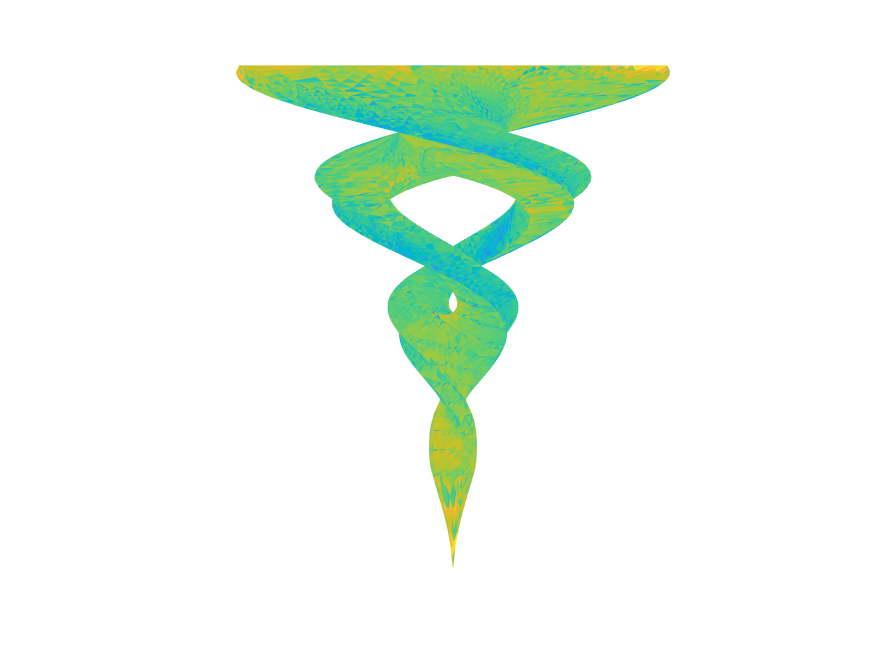} &
        \includegraphics[clip,trim = {4cm 1cm 3.8cm 0.8cm},width = 0.15\textwidth]{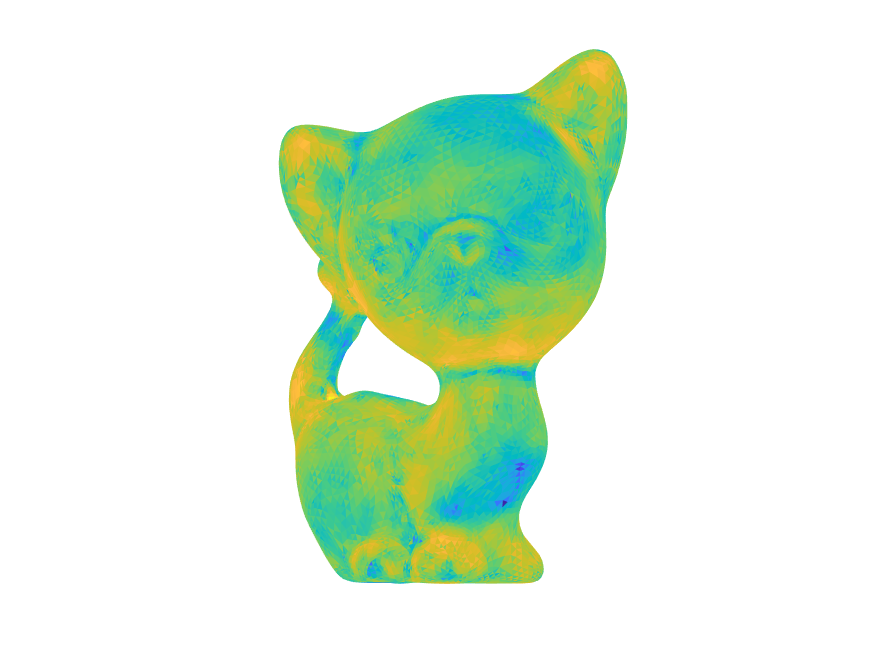} &
        \includegraphics[clip,trim = {4cm 1cm 3.8cm 0.8cm},width = 0.15\textwidth]{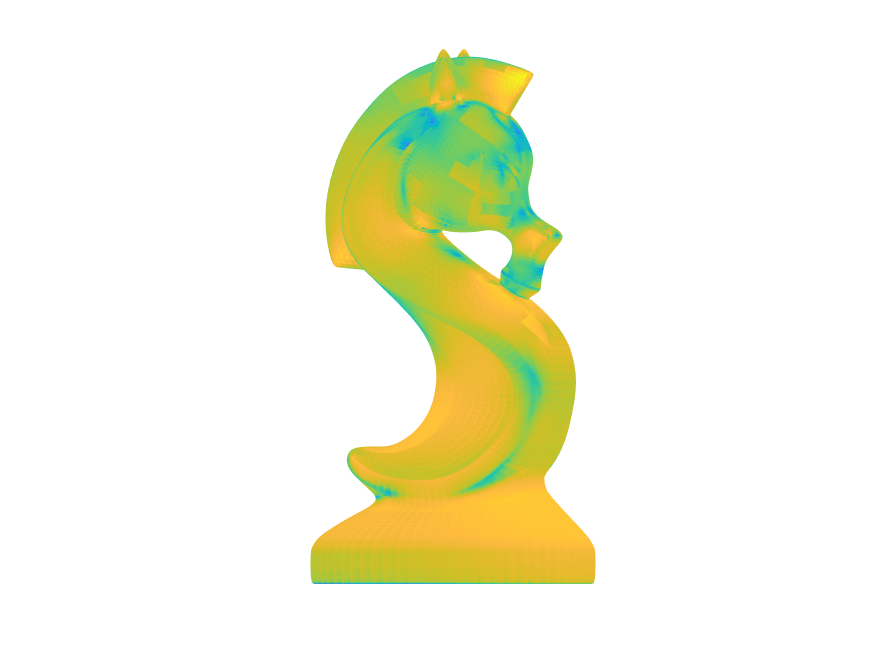} &
        \includegraphics[clip,trim = {4cm 2cm 3.5cm 0.8cm},width = 0.15\textwidth]{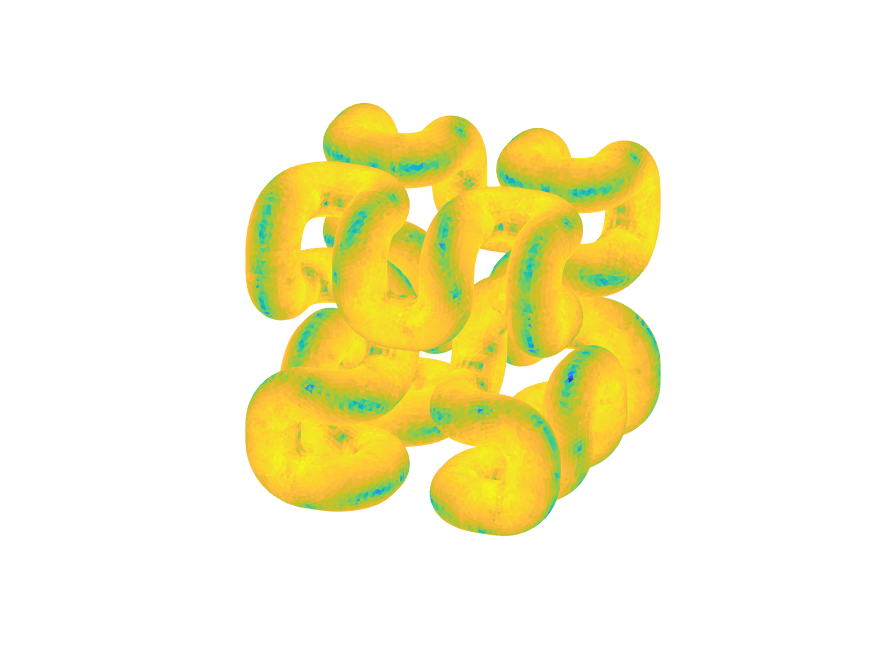} &
        \includegraphics[clip,trim = {3.cm 1cm 2.8cm 1cm},width = 0.15\textwidth]{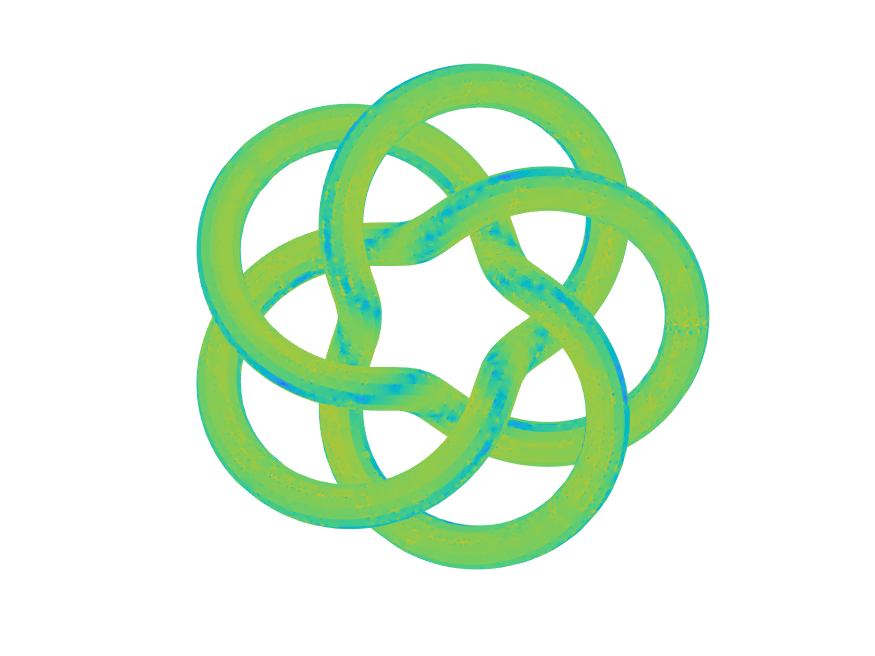} &
        \includegraphics[clip,trim = {5cm 2cm 3.8cm 2cm},width = 0.15\textwidth]{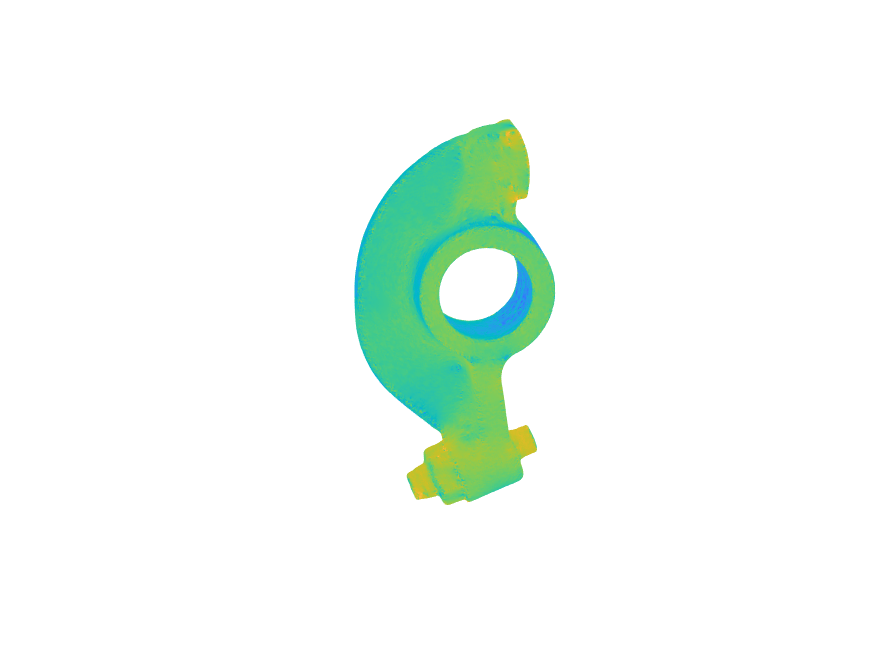} &
        \includegraphics[clip,trim = {12cm 6cm 0cm 0cm}, width = 0.1\textwidth]{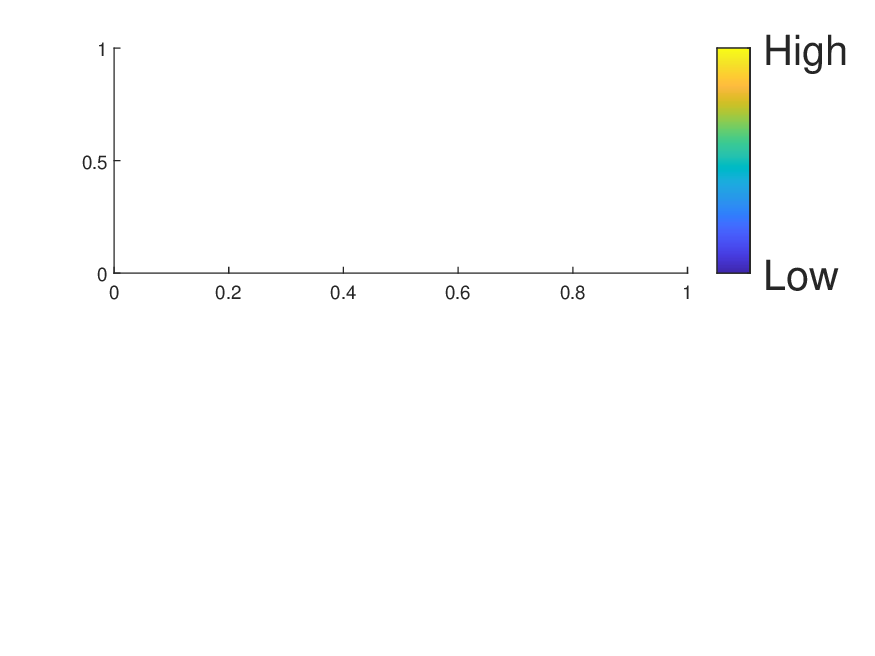} \\
        Twirl & Kitten & ChessHorse & Hilb64Thick & Knot & RockerArm &
    \end{tabular}}
    \caption{Genus-one mesh models with handle and tunnel loops computed by ReebHanTun \cite{TDFF13} and the angle distortion distributions of the maps resulting from DPCF. The magnitude of the value is relative to the maximum distortions on the surfaces.}
    \label{fig:g1_ori}
\end{figure}

\begin{figure}[htp]
    \centering
    \resizebox{\textwidth}{!}{
    \begin{tabular}{cccccc}
        \includegraphics[clip,trim = {0cm 0cm 0cm 0cm},width = 0.3\textwidth]{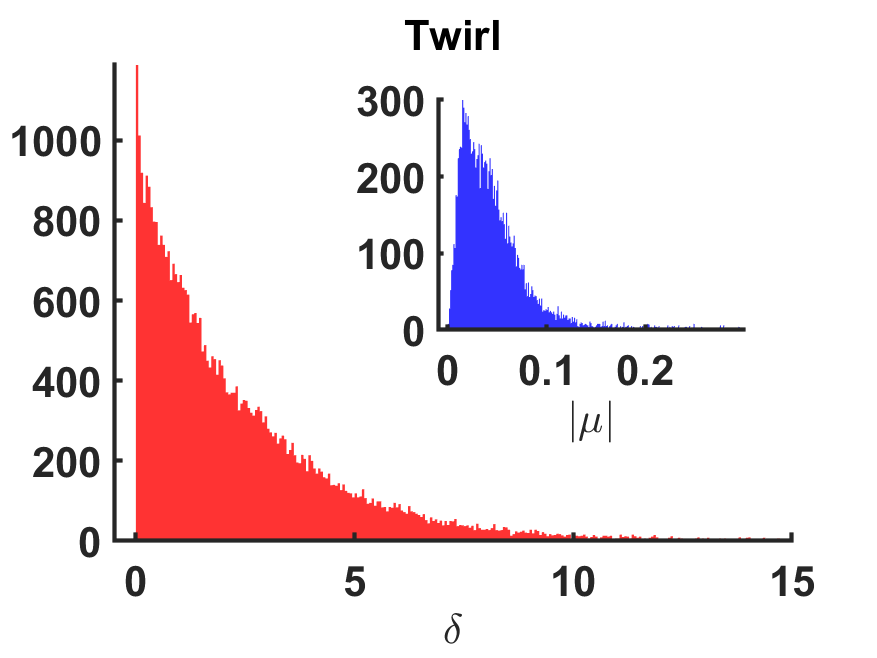} &
        \includegraphics[clip,trim = {0cm 0cm 0cm 0cm},width = 0.3\textwidth]{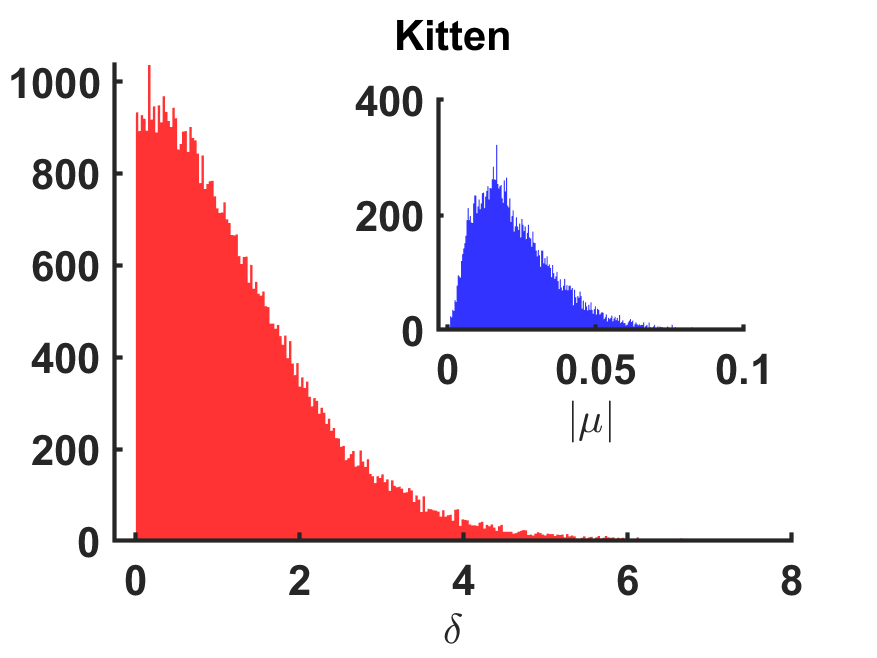} &
        \includegraphics[clip,trim = {0cm 0cm 0cm 0cm},width = 0.3\textwidth]{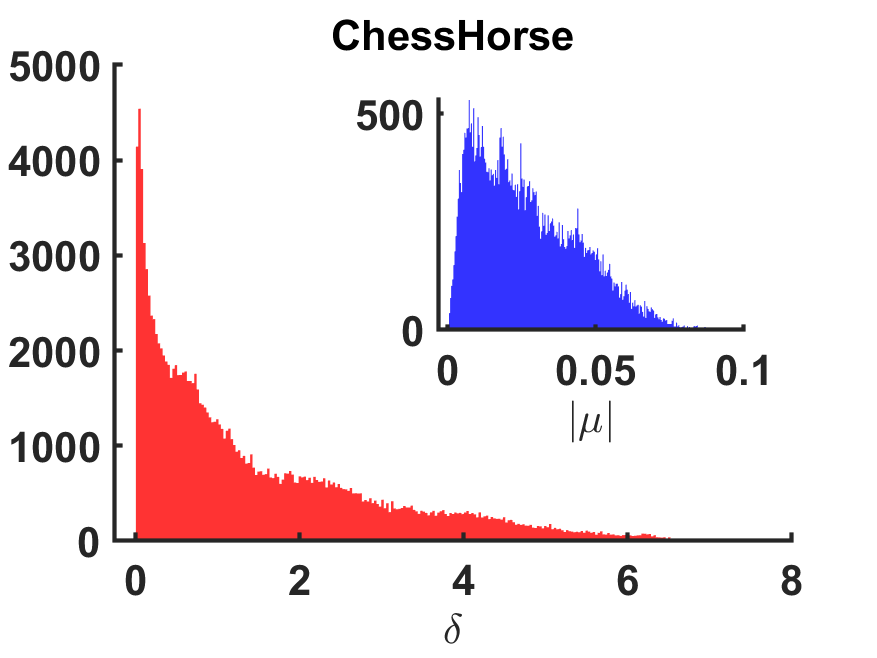} \\
        \includegraphics[clip,trim = {0cm 0cm 0cm 0cm},width = 0.3\textwidth]{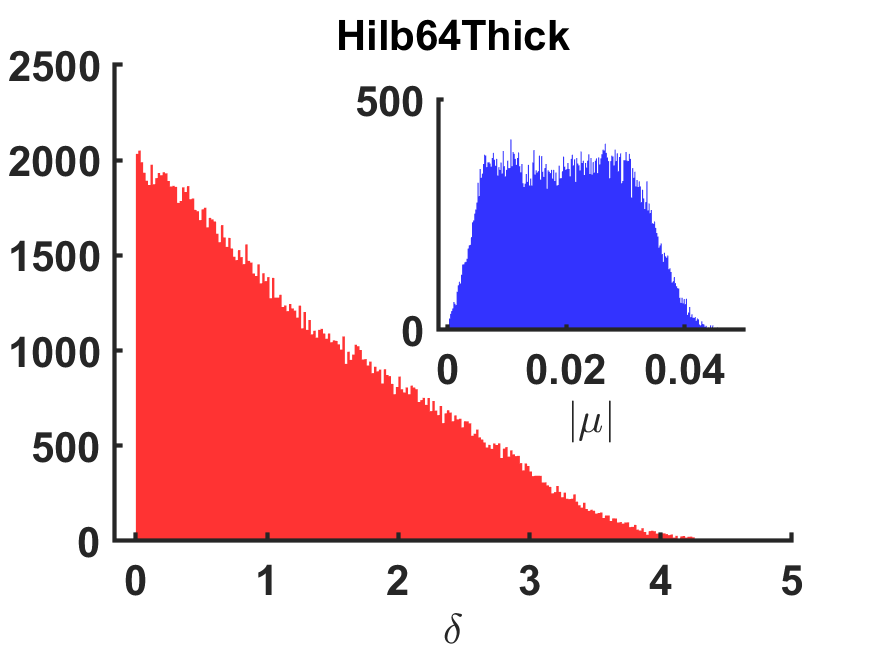} &
        \includegraphics[clip,trim = {0cm 0cm 0cm 0cm},width = 0.3\textwidth]{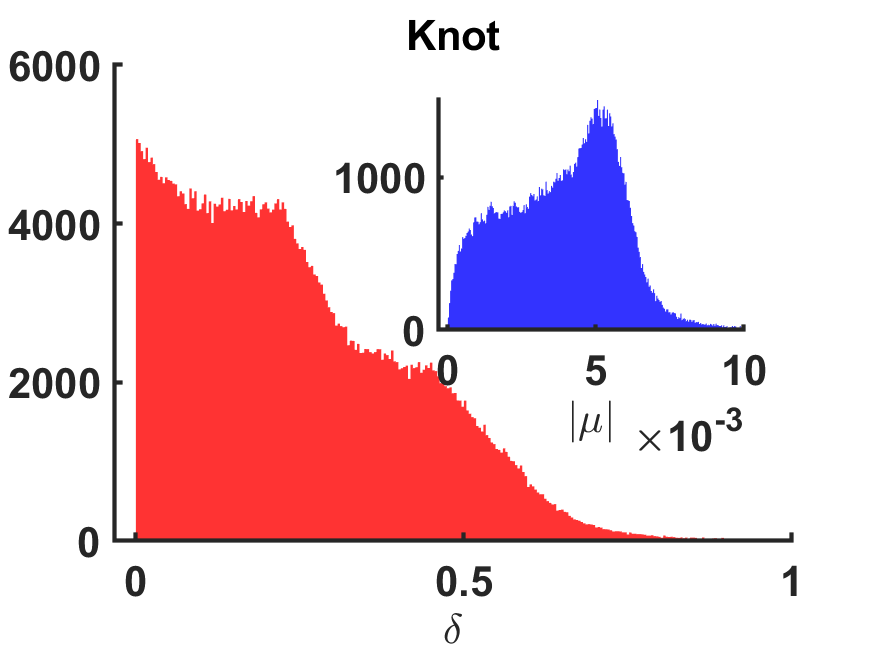} &
        \includegraphics[clip,trim = {0cm 0cm 0cm 0cm},width = 0.3\textwidth]{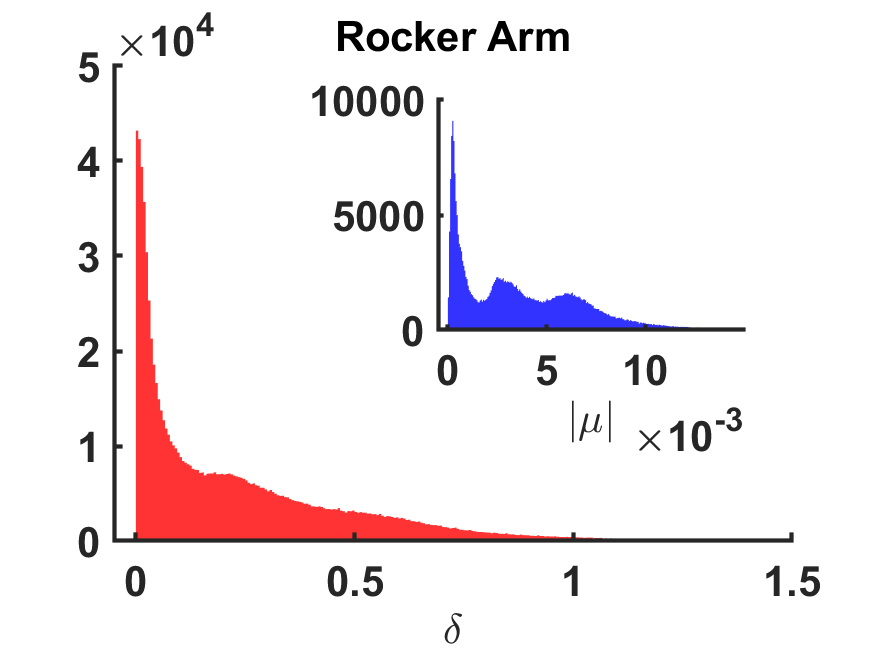} &  \\
    \end{tabular}}
    \caption{Angle distortion histograms of the maps resulting from DPCF. The red figures represent absolute angle errors (degree), while the blue ones represent Beltrami coefficients.}
    \label{fig:g1_hist}s
\end{figure}

The genus-one surface meshes used for the experiment and their basic information are shown in the first rows of Figure \ref{fig:g1_ori} and Table \ref{tab:g1_ori}. The red and blue curves displayed in Figure \ref{fig:g1_ori} represent the handle and tunnel loops computed by ReebHanTun \cite{TDFF13}, respectively. In Table \ref{tab:g1_ori}, $\#$ face and $\#$ vert. represent the numbers of triangle faces and vertices, respectively.
We compare our results with those of Algorithm 1 in \cite{MHTL20}, which is based on de Rham cohomology theory.
Surprisingly, we test more than 20 meshes and find that their conformal accuracies are mostly identical.
The relative error of DPCF is $0.2\%$ slightly better than that of Algorithm 1 in \cite{MHTL20}.
Hence, we only show the conformal distortions, i.e., $\delta$ and $|\mu|$ as in \eqref{err:aae} and \eqref{err:Beltrami}, of DPCF in Table \ref{tab:g1_ori} and Figure \ref{fig:g1_hist}.
One can see that the conformal distortions are small and most of them are close to 0.
However, DPCF is almost $5$ times more efficient than the competing method.
Additionally, Figure \ref{fig:g1_ori} shows the angle distortion $\delta$ distributions of DPCF on 6 meshes, in which all the quantities in degree are quite small according to Table \ref{tab:g1_ori}.
One can see that no additional conformal distortions occur near the cutting paths, which is independent of the cutting paths of DPCF.

\begin{table}[htp]
    \centering
    \resizebox{\textwidth}{!}{
    \begin{tabular}{c|cc|cc|cc|cc|c}
    \hline
        \multirow{2}{*}{Mesh} & \multirow{2}{*}{$\#$ face} & \multirow{2}{*}{$\#$ vert.} & \multicolumn{2}{c|}{$\delta$ (degree)} & \multicolumn{2}{c|}{$100*|\mu|$} & \multicolumn{3}{c}{Time (s)} \\ \cline{4-10}
        &&& Mean & Std & Mean & Std & alg.1\cite{MHTL20} & DPCF & ratio \\
    \hline
        Twirl & 14208 & 7104 & 2.762 & 4.697 & 5.688 & 8.103 & 0.147 & \textbf{0.025} & 5.853 \\
        Kitten & 20000 & 10000 & 1.270 & 1.080 & 2.317 & 1.459 & 0.210 & \textbf{0.034} & 6.091 \\
        ChessHorse & 46016 & 23008 & 1.521 & 1.484 & 2.760 & 1.924 & 0.527 & \textbf{0.085} & 6.173 \\
        Hilb64Thick & 64044 & 32022 & 1.249 & 0.929 & 2.000 & 1.000 & 0.946 & \textbf{0.153} & 6.174 \\
        Knot & 169532 & 84766 & 0.248 & 0.230 & 0.393 & 0.298 & 3.074 & \textbf{0.553} & 5.563 \\
        RockerArm & 309646 & 154823 & 0.275 & 0.385 & 0.438 & 0.517 & 5.332 & \textbf{0.965} & 5.528 \\
    \hline
    \end{tabular}
    }
    \caption{Comparison between Algorithm 1 in \cite{MHTL20} and DPCF on genus-one surfaces. $\#$ face and $\#$ vert. represent the number of triangle faces and vertices, respectively. $\delta$ is the absolute angle error and $\mu$ is the Beltrami coefficient defined in \eqref{err:aae} and \eqref{err:Beltrami}, respectively. }
    \label{tab:g1_ori}
\end{table}

\subsection{Bijectivity Guarantee Strategy}

As in Theorem \ref{thm:delG1}, DPCF guarantees the bijectivity provided that the mesh $\mathcal{M}$ is intrinsic Delaunay, i.e., the sum of the opposite angles in the adjacent triangles is no more than $\pi$ for each interior edge, which means $w_{ij} > 0$.
Furthermore, \cite{ABBS07} demonstrated that a non-Delaunay surface triangulation can be transformed into an intrinsic Delaunay triangulation (IDT) through iterative edge flipping.
Hence, a direct approach for guaranteeing the bijectivity of the resulting map is to modify the mesh $\mathcal{M}$ to be the associated IDT $\mathcal{M}_d$ by executing edge flipping first and then performing conformal flattening.
Since the triangular mesh is the discrete representation of the corresponding continuous surface,
the resulting maps on the original mesh $\mathcal{M}$ and the associated IDT $\mathcal{M}_d$ generally converge to an identical continuous conformal map.
To observe the influence of the IDT preprocessing strategy, we show the performance achieved by DPCF without and with IDT preprocessing on \emph{Torus} in Table \ref{tab:IDT}.
This model is remeshed to two resolutions with low-quality so that there exist folding triangles on flattened domains.
We can see that the IDT preprocessing strategy can eliminate all folding triangles and provide better results with lower conformal distortions.

\begin{table}[htp]
\centering
\begin{tabular}{c|cc|c|cc|cc|c}
    \hline
        \multirow{2}{*}{Mesh} & \multirow{2}{*}{$\#$ face} & \multirow{2}{*}{$\#$ vert.} & \multirow{2}{*}{IDT} & \multicolumn{2}{c|}{$\delta$ (degree)} & \multicolumn{2}{c|}{$100*|\mu|$} & \multirow{2}{*}{$\#$ folding} \\ \cline{5-8}
        &&&& Mean & Std & Mean & Std & \\
    \hline
\multirow{4}{*}{Torus} & \multirow{2}{*}{16920} & \multirow{2}{*}{8460} & \ding{55} & 0.726 & 1.267 & 1.663 & 4.464 & 10 \\
 & & & \ding{51} & \textbf{0.479} & \textbf{0.516} & \textbf{0.820} & \textbf{0.766} & 0 \\
 \cline{2-9}
& \multirow{2}{*}{32500} & \multirow{2}{*}{16250} & \ding{55} & 0.626 & 1.289 & 1.510 & 4.754 & 30 \\
& & & \ding{51} & \textbf{0.351} & \textbf{0.427} & \textbf{0.600} & \textbf{0.831} & 0 \\
\hline
\end{tabular}
\caption{Comparison between DPCF without and with IDT preprocess on \emph{Torus}.}
\label{tab:IDT}
\end{table}

\section{Applications} \label{sec:app}

\subsection{Texture Mapping}

\begin{figure}[htp]
    \centering
    \resizebox{0.8\textwidth}{!}{
    \begin{tabular}{cccc}
        \includegraphics[clip,trim = {4cm 1cm 3.8cm 0.8cm},width = 0.23\textwidth]{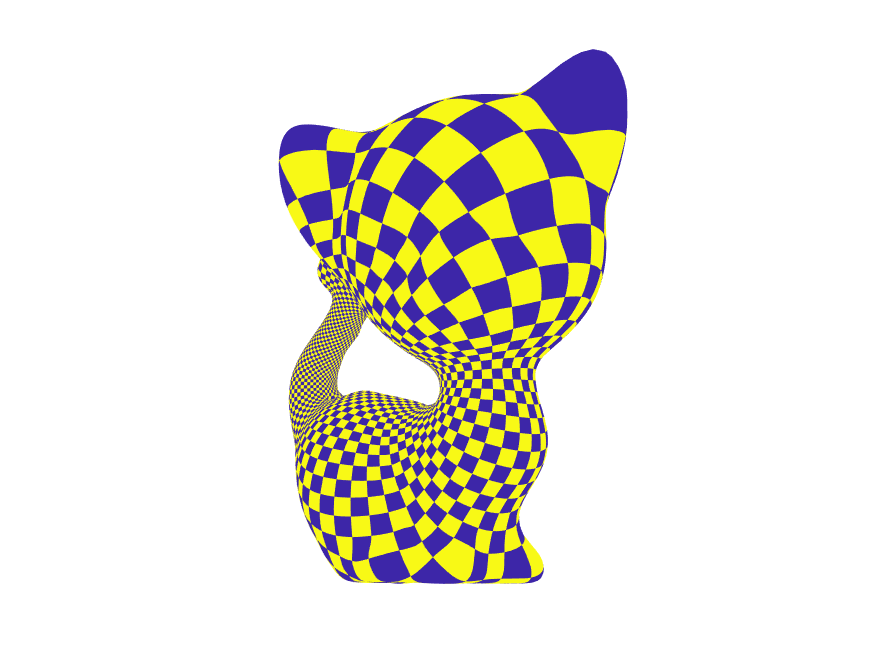} &
        \includegraphics[clip,trim = {5cm 2cm 3.8cm 2cm},width = 0.23\textwidth]{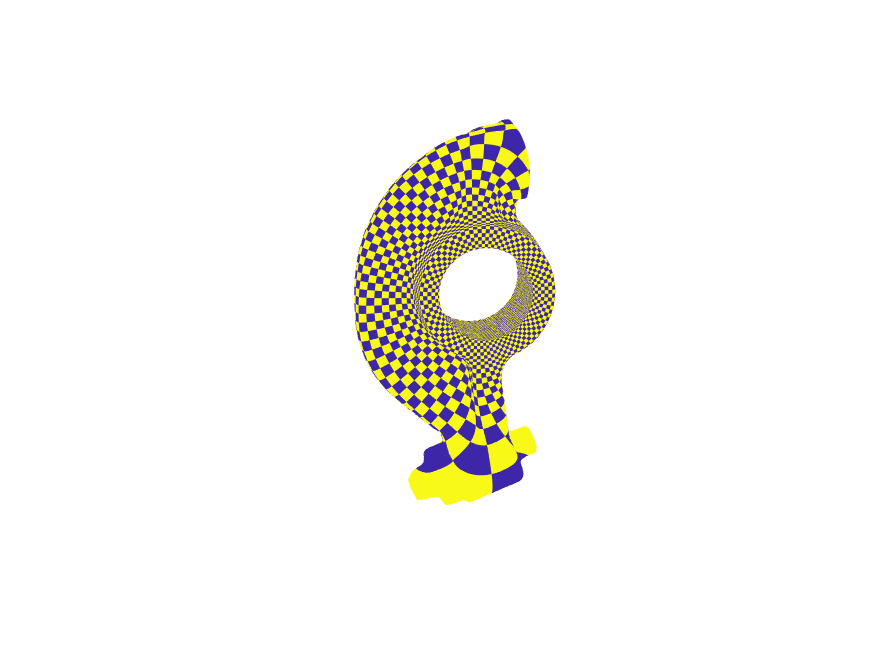} &
        \includegraphics[clip,trim = {2cm 1cm 1.8cm 0.8cm},width = 0.32\textwidth]{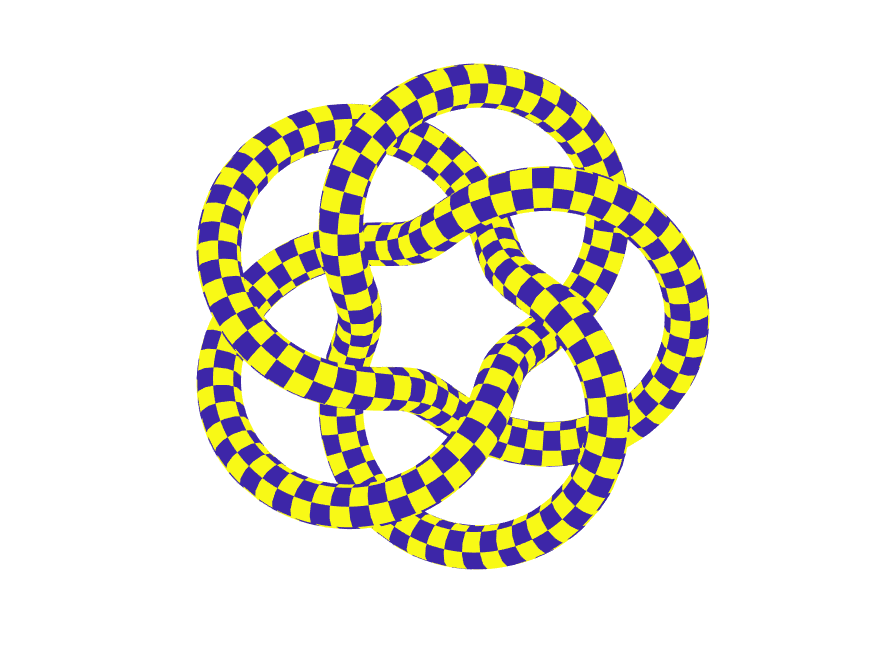}
        \\
        Kitten & RockerArm & Knot
    \end{tabular}
    }
    \caption{Chessboard texture mappings on the genus-one surfaces, \emph{Kitten}, \emph{RockerArm} and \emph{Knot}.}
    \label{fig:texture}
\end{figure}

A common method for achieving texture mapping is to flatten the target surface first and then design the desired texture on the flattened domain.
Finally, the texture appears on the surface through the inverse flattening map.
The main advantage of applying a conformal map on texture mapping is its local shape preservation property, so designers can intuitively observe the effect of their texture on the plane.
Moreover, owing to the seamless translation correspondence between the cutting paths, we do not have to account for the additional need to seamlessly glue the cutting path edges.
Figure \ref{fig:texture} depicts chessboard texture mappings on genus-one surfaces. The textures near the cutting paths are well rendered and the orthogonality of the curves is preserved.

\subsection{Vertebra Remeshing}

While the initial bone and organ surface meshes extracted from CT segmentation captures the overall vertebral geometry, it typically suffers from several deficiencies that compromise its suitability for downstream numerical simulations.
Poorly shaped and highly non-uniform vertex distributions are common artefacts that can lead to numerical instability or excessive computational cost in medical analysis.
To obtain a robust and physiologically meaningful biomechanical response, remeshing the surface is therefore an essential pre-processing step, aimed at generating a high-quality and curvature-adaptive mesh while preserving critical anatomical features.

\begin{figure}
    \centering
    \includegraphics[width=0.9\linewidth]{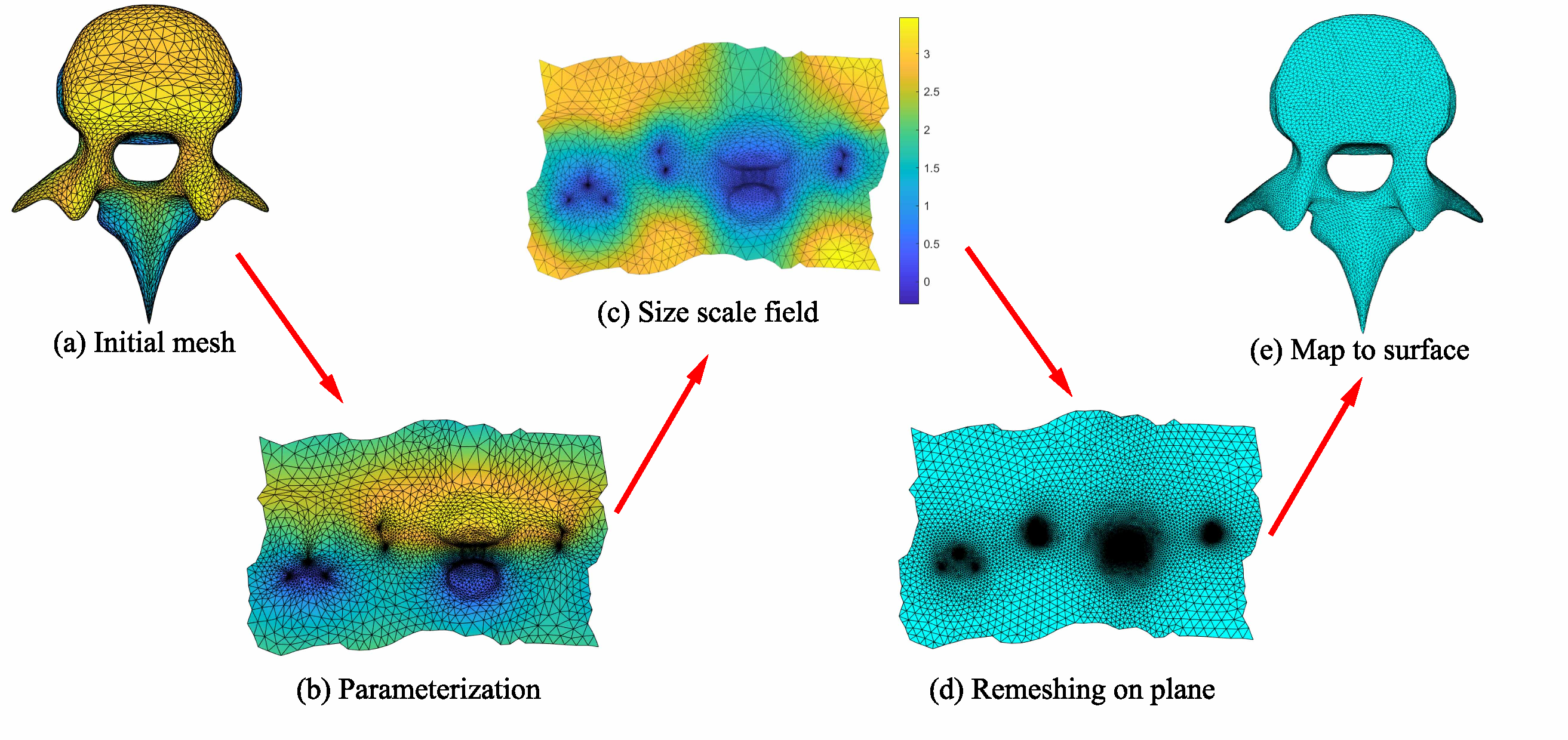}
    \caption{Vertebra remeshing proccess.}
    \label{fig:VertebraRemesh}
\end{figure}

Here we take the vertebra, a genus-one surface, as an example and perform remeshing via our proposed DPCF algorithm.
Given an initial vertebra mesh $\mathcal{M}$ as in Figure \ref{fig:VertebraRemesh} (a), which has 7,500 faces and 3,750 vertices, we can flatten it comformally by DPCF, as in Figure \ref{fig:VertebraRemesh} (b). This conformal map induces a size scale field (Figure \ref{fig:VertebraRemesh} (c)), which measures the local size distortion from the conformal map. Through the size scale field, we can determine the local size of the prescribed triangle element. Then the adaptive mesh generation with boundary correspondence is performed accordingly, as in Figure \ref{fig:VertebraRemesh} (d).
Finally, by using inverse conformal map, the reshemed vertebra surface is obtained, as in Figure \ref{fig:VertebraRemesh} (e), which has 28,192 faces and 14,096 vertices.

\begin{figure}
    \centering
    \includegraphics[width=0.4\linewidth]{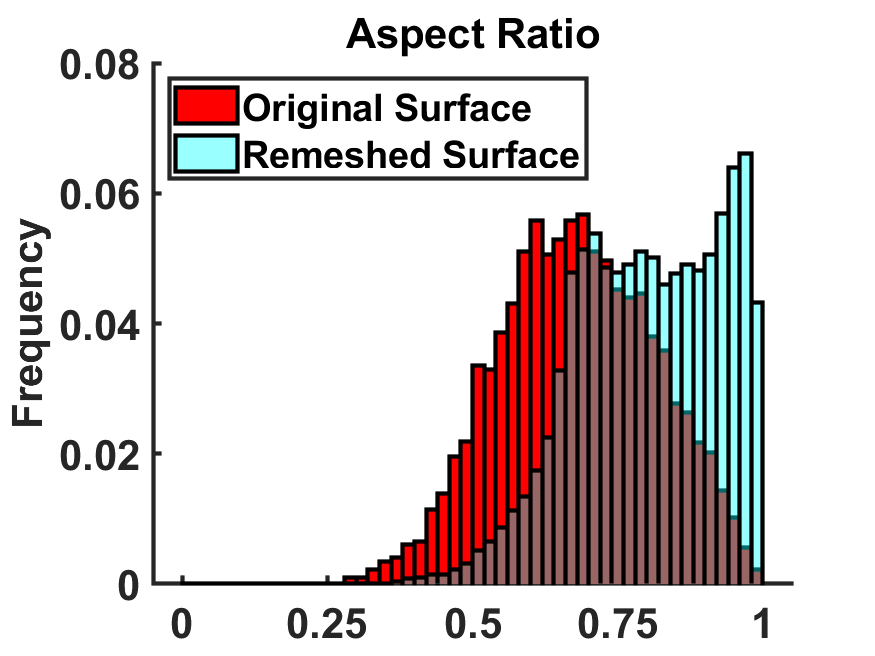}
    \includegraphics[width=0.4\linewidth]{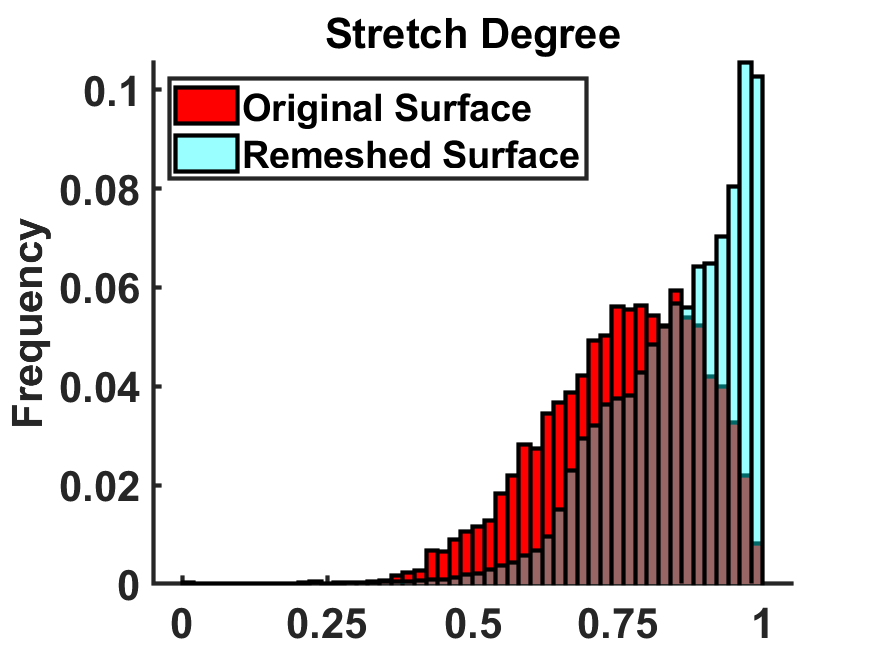}
    \caption{Aspect ratio and stretch degree histograms of initial and remeshed surface meshes.}
    \label{tab:QemeshQuality}
\end{figure}

Thanks to the angle-preservation property of conformal map, the high-quality planar triangular elements generated in the parameter domain are largely preserved when mapped back to the vertebral surface.
Consequently, the remeshed surface also consists of high-quality elements.
Figure \ref{tab:QemeshQuality} illustrates the mesh quality histograms before and after remeshing, where x-axis represents the mesh quality value and y-axis represents the frequency.
The mesh quality is measured by aspect ratio $Q_\text{AR}$ and stretch degree $Q_\text{SD}$, which are defined as
\begin{align}
    Q_\text{AR}(\tau) = \frac{\ell_\text{min}}{\ell_\text{max}}, \quad
    Q_\text{SD}(\tau) = \frac{\sqrt{12}r}{\ell_\text{max}},
\end{align}
respectively, with $\ell_\text{min}$, $\ell_\text{max}$ denoting the minimal and maximal edge lengths of triangle $\tau$, respectively, and $r$ being its inradius.
Both $Q_\text{AR}$ and $Q_\text{SD}$ lie in the interval $[0,1]$ and the closer they are to $1$, the better the mesh quality.
The histograms clearly demonstrate a substantial improvement in mesh quality after remeshing, with the vast majority of elements exhibiting $Q_\text{AR}$ and $Q_\text{SD}$ close to $1$.

\section{Conclusion}\label{sec:conclusion}

In this paper, we propose a novel periodic conformal flattening method for genus-one surfaces via conformal energy minimization with seam constraints, which can be expressed as a quadratic functional minimization problem with respect to the vertices and translation vectors.
Naturally, the conformal map associated with the optimal boundary is obtained by solving a linear system, which leads to the high efficiency of our proposed DPCF algorithm.
Theoretically, we prove that the maps resulting from our proposed DPCF do not depend on the cutting path and that no additional error occurs near the cutting path.
Furthermore, we can prove that the translation vectors rely only on the homology basis to which the cutting paths belong, which illustrates the stability of the proposed algorithms in cutting path selection scenarios.
Additionally, we prove that the resulting map is bijective under the positive edge weight condition, which inspires an intrinsic Delaunay strategy for guaranteeing bijectivity.
In numerical experiments, we present the performance of DPCF on genus-one surfaces. The comparison with the de Rham cohomology approach \cite{MHTL20} illustrates the low conformal distortion of our approach and at least 5-time efficiency improvements, respectively. The angle distortion distributions demonstrate the seamless property of our proposed algorithms.
For general surfaces, the intrinsic Delaunay preprocess can guarantee the bijectivety and high accuracy, which can be applied to low-quality meshes.
On the basis of this strategy, we illustrate the applications on texture mapping and medical imaging to demonstrate the practicality of our method.

In our method, conformal energy minimization plays an important role and shows its effectiveness in surface conformal flattening, especially in terms of efficiency.
The periodic flattening method proposed herein is applicable for surfaces that can be flattened to the fundamental domains.
Therefore, this approach cannot be directly applied to arbitrary surfaces, such as genus-zero and higher-genus surfaces, which are flattened with more complicated seam constraints.
The development of efficient seamless conformal flattening algorithms with greater generality through conformal energy minimization is one of our ongoing research goals.

\bibliographystyle{plain}

\end{document}